\documentstyle[12pt]{article}
\newtheorem{theorem}{Theorem}[section]
\newtheorem{lemma}[theorem]{Lemma}
\newtheorem{proposition}[theorem]{Proposition}

\newtheorem{corollary}[theorem]{Corollary}
\newtheorem{remar}[theorem]{Remark}
\newenvironment{proof}{Proof:\ \ \ }{\QED}
\newenvironment{remark}{\begin{remar}\rm}{\end{remar}}

\newcommand{\QED}{{\unskip\nobreak\hfil\penalty50%
\hskip1em\hbox{}\nobreak\hfil $\Box$%
\parfillskip=0pt \finalhyphendemerits=0 \par\medskip\noindent}}
\newcommand{\bfind}[1]{\index{#1}{\bf #1}}

\newcommand{\n}{\par\noindent}

\newcommand{\sn}{\par\smallskip\noindent}

\newcommand{\bn}{\par\bigskip\noindent}
\newcommand{\pars}{\par\smallskip}
\newcommand{\parm}{\par\medskip}
\newcommand{\parb}{\par\bigskip}

\newcommand{\ovl}[1]{\overline{#1}}

\newcommand{\sep}{^{\rm sep}}
\newcommand{\iT}{^{[i]}}
\newcommand{\chara}{\mbox{\rm char}\,}
\newcommand{\trdeg}{\mbox{\rm trdeg}\,}
\newcommand{\subsetuneq}{\mathrel{\raisebox{.8ex}{\footnotesize%
$\displaystyle\mathop{\subset}_{\not=}$}}}
\newcommand{\adresse}{\par\bigskip \small\rm
 Department of Mathematics and Statistics, 
 University of Saskatchewan, \par
 106 Wiggins Road, 
 Saskatoon, Saskatchewan, Canada S7N 5E6 \par
 email: fvk@math.usask.ca \ \ --- \ \ home page:
http://math.usask.ca/$\,\tilde{ }\,$fvk/index.html}
\font\tenlv=msbm10 scaled 1200
\font\sevenlv=msbm7 scaled 1200
\font\fivelv=msbm5 scaled 1200
\def\lv #1{{\mathchoice{{\hbox{\tenlv #1}}}{{\hbox{\tenlv #1}}}
{{\hbox{\sevenlv #1}}}{{\hbox{\fivelv #1}}}}}

\newcommand{\N}{\lv N}
\newcommand{\Q}{\lv Q}
\newcommand{\R}{\lv R}
\newcommand{\Z}{\lv Z}

\begin{document}
\title{On local uniformization in arbitrary characteristic I} 
\author{Franz--Viktor Kuhlmann}
\date{23.\ 10.\ 1998}
\maketitle
\begin{abstract}\noindent
We prove that every place of an algebraic function field $F|K$ of
arbitrary characteristic admits local uniformization in a finite
extension ${\cal F}$ of $F$. We show that ${\cal F}|F$ can be chosen
to be normal. If $K$ is perfect and $P$ is of rank 1, then
alternatively, ${\cal F}$ can be obtained from $F$ by at most two Galois
extensions; if in addition $P$ is zero-dimensional, then we only need
one Galois extension. Certain rational places of rank 1 can be
uniformized already on $F$. We introduce the notion of ``relative
uniformization'' for arbitrary finitely generated extensions of valued
fields. Our proofs are based solely on valuation theoretical theorems,
which are of fundamental importance in positive characteristic.
\end{abstract}
%
%
%
%
\section{Introduction}
In [Z1], Zariski proved the Local Uniformization Theorem for places of
algebraic function fields over base fields of characteristic 0. In [Z3],
he uses this theorem to prove resolution of singularities for surfaces
in characteristic 0. As the resolution of singularities for arbitrary
dimensions in positive characteristic is still an open problem, one is
interested in generalizations of the Local Uniformization Theorem
to positive characteristic (cf.\ [S]). In this paper, we will prove a
weak version of the Local Uniformization Theorem, for function fields of
arbitrary characteristic:
\begin{theorem}                                \label{MT}
Let $F|K$ be a function field of arbitrary characteristic,
and $P$ a place of $F|K$. Then there exist a finite normal extension
${\cal F}$ of $F$, an extension of $P$ from $F$ to ${\cal F}$,
a finite purely inseparable extension ${\cal K}$ of
$K$ within ${\cal F}$ and a model of ${\cal F}|{\cal K}$ on
which $P$ is centered in a smooth point.
\end{theorem}
Throughout this paper, ``function field'' will always mean ``algebraic
function field''. By a place of $F|K$ we mean a place whose restriction
to $K$ is the identity. Talking of an extension ${\cal F}$ of $F$, we
will from now on tacitly assume that it is equipped with an extension of
$P$ (which is again denoted by $P$).

Theorem~\ref{MT} follows from the results of J.~de Jong [dJ] (who proves
resolution of singularities after a finite normal extension of the
function field). However, we will give an entirely valuation theoretical
proof which will provide important additional information. In
particular, we wish to get as close as possible to taking ${\cal F}|F$
Galois. Why do we want that ${\cal F}|F$ is Galois? Apart from
geometrical reasons, it is because the work of S.~Abhyankar seems to
indicate that there is a chance to ``pull down'' local uniformization
through Galois extensions. This would give us what we actually want:
local uniformization without extending the function field.

On the other hand, for certain applications of Theorem~\ref{MT} (e.g.,
to the model theory of fields in the spirit of [J--R]; cf.\ also [K3]),
it is important to have a valuation theoretical control on the extension
${\cal F}|F$ and the residue field extension ${\cal F}P|FP$. (We want to
have ${\cal F}P$ to be as close to $FP$ as possible, but in positive
characteristic we may expect that we have to take a purely inseparable
extension into the bargain.) We cannot obtain this control if we insist
that ${\cal F}|F$ be Galois. Instead, we will show in a subsequent
paper [K6] that in case of a perfect base field $K$, the extension
${\cal F}|F$ can always be chosen to be separable and such that
${\cal F} P|FP$ is purely inseparable. See also [K7] for background
information.

\parm
In order to prove Theorem~\ref{MT}, we will give a description of the
special form in which the Jacobian condition for smoothness can be
satisfied. For polynomials $f_1,\ldots,f_n$ in variables $X_1,\ldots,
X_n\,$, we write $f=(f_1,\ldots,f_n)$ and denote by $J_f$ the Jacobian
matrix $\left(\frac{\partial f_i}{\partial X_j}\right)_{i,j}$ of their
partial derivatives. Take a finitely generated extension ${\cal F}|
{\cal K}$, not necessarily transcendental, and a place $P$ of
${\cal F}$, not necessarily the identity on ${\cal K}$. We write
${\cal O}_{\cal F}$ for the valuation ring of $P$ on ${\cal F}$, and
${\cal O}_{\cal K}$ for that of $P$ on ${\cal K}$. For an element
$a$, its $P$-residue will be denoted by $aP$. We will say that
\bfind{$({\cal F}|{\cal K},P)$ is weakly uniformizable} if there are
\pars
$\bullet$ \ a transcendence basis $T=\{t_1,\ldots,t_s\}\subset
{\cal O}_{\cal F}$ of ${\cal F}|{\cal K}\,$ (which may be empty),\par
$\bullet$ \ elements $\eta_1,\ldots,\eta_n\in {\cal O}_{\cal F}\,$,\par
$\bullet$ \ polynomials $f_i(X_1,\ldots,X_n)\in {\cal O}_{\cal K}
[t_1,\ldots,t_s,X_1,\ldots,X_n]$, $1\leq i\leq n$,
\sn
such that ${\cal F}={\cal K}(t_1,\ldots,t_s,\eta_1,\ldots,\eta_n)$,
and\pars
(U1) \ for $i<j$, $X_j$ does not occur in $f_i\,$,\par
(U2) \ $f_i(\eta_1,\ldots,\eta_n)=0$ for $1\leq i\leq n$,\par
(U3) \ $(\det J_f (\eta_1,\ldots,\eta_n))P\ne 0$.
\sn
In this case, we will call $T$ a \bfind{uniformizing transcendence
basis}.

\pars
Assertion (U1) implies that $J_f$ is lower triangular. Assertion (U3)
says that
\begin{equation}                            \label{Jac}
\det J_{fP} (\eta_1P,\ldots,\eta_nP)\>=\>
(\det J_f (\eta_1,\ldots,\eta_n))P \>\ne\> 0\;,
\end{equation}
where $fP=(f_1P,\ldots,f_nP)$ and $f_iP$ denotes the $P$-reduction of
$f_i\,$, i.e., the polynomial obtained from $f_i$ through replacing
every coefficient by its $P$-residue. Note that $\det J_f (\eta_1,
\ldots,\eta_n)\in {\cal O}_{\cal F}$ since $t_1,\ldots,t_s,
\eta_1,\ldots,\eta_n\in {\cal O}_{\cal F}$ and the $f_i$ have
coefficients in ${\cal O}_{\cal K}\,$.

Given elements $\zeta_1,\ldots,\zeta_m\in {\cal O}_{\cal F}$, we will
say that \bfind{$({\cal F}|{\cal K},P)$ is uniformizable with respect
to $\zeta_1,\ldots,\zeta_m$} if the elements $\eta_1,\ldots,\eta_n$ can
be chosen such that $\zeta_1,\ldots,\zeta_m$ appear among them. We
say that \bfind{$({\cal F}|{\cal K},P)$ is uniformizable} if it is
uniformizable with respect to {\it every} choice of
finitely many elements in ${\cal O}_{\cal F}\,$. This property is
transitive:
\begin{theorem}                             \label{translu}
If $({\cal F}|{\cal L},P)$ and $({\cal L}|{\cal K},P)$ are
uniformizable, then $({\cal F}|{\cal K},P)$ is uniformizable.
\end{theorem}
Note the following well-known fact: if all conditions except (U1) for
weak uniformizability hold for $g_i$ in the place of $f_i\,$, then there
are $f_i$ which satisfy all conditions (for the same $T$ and $\eta$'s).
We only have included condition (U1) since it is a nice natural
side-effect of our approach which uses the transitivity.

\pars
Now assume that ${\cal F}|{\cal K}$ is a function field (i.e., $\trdeg
{\cal F}|{\cal K} \geq 1$) and that $P$ is the identity on ${\cal K}$.
Then ${\cal O}_{\cal K}= {\cal K}$, and the $P$-residues of the
coefficients are obtained by just replacing $t_j$ by $t_jP$, for $1\leq
j\leq n$. Hence if we view the polynomials $f_i$ as polynomials in the
variables $Z_1,\ldots,Z_s, X_1,\ldots,X_n$ with coefficients in
${\cal K}$, then (\ref{Jac}) means that at the point
$(t_1P,\ldots,t_sP,\eta_1P,\ldots,\eta_nP)$ the diagonal elements of
$J_f$ and thus also its determinant do not vanish. This assertion says
that on the variety defined over ${\cal K}$ by the $f_i$ (and having
generic point $(t_1,\ldots,t_s,\eta_1,\ldots, \eta_n)$ and function
field ${\cal F}$), the place $P$ is centered at the smooth point
$(t_1P,\ldots,t_sP,\eta_1P,\ldots, \eta_nP)$.

\pars
This discussion shows that Theorem~\ref{MT}
is a consequence of the following two theorems:
\begin{theorem}                                \label{MT1}
Let $F|K$ be a function field of arbitrary characteristic
and $P$ a place of $F|K$. Take any
elements $\zeta_1,\ldots,\zeta_m$ in the valuation ring ${\cal O}_F$ of
$P$ on $F$. Then there exist a finite extension ${\cal F}$ of $F$, an
extension of $P$ to ${\cal F}$, and a finite purely inseparable
extension ${\cal K}$ of $K$ within ${\cal F}$ such that
$({\cal F}|{\cal K},P)$ is uniformizable with respect to
$\zeta_1,\ldots,\zeta_m\,$.
%
%
\end{theorem}
(Here, ``$F.{\cal K}$'' denotes the field compositum of ${\cal K}$ and
$F$ inside of ${\cal F}$, i.e., the smallest subfield of ${\cal F}$
containing ${\cal K}$ and $F$.) Note that ${\cal K}=K$ if $K$ is
perfect.
\begin{theorem}                             \label{MT1norm}
Take any subextension $E|K$ of $F|K$ of the same transcendence degree.
Then in addition to the assertion of Theorem~\ref{MT1}, ${\cal F}$ can
always be chosen to be a normal extension of $E$ and of $F$.
\end{theorem}

\pars
By uniformizing with respect to the $\zeta$'s, we obtain the following
important information: if we have already a model $V$ of $F|K$ with
generic point $(z_1,\ldots,z_k)$, where $z_1,\ldots,z_k\in {\cal O}
_F\,$, then we can choose our new model ${\cal V}$ of ${\cal F}|
{\cal K}$ in such a way that the local ring of the center of $P$ on
${\cal V}$ contains the local ring of the center $(z_1P,\ldots,z_kP)$ of
$P$ on $V$. For this, we only have to let $z_1,\ldots,z_k$ appear among
the $\zeta$'s.

\pars
In important special cases, we can show much stronger results. Before
we state them, let us introduce some useful notions. Let $P$ be an
arbitrary place on a field $F$. We will call $(F,P)$ a \bfind{valued
field}, keeping in mind its associated valuation, which we denote by
$v_P^{ }\,$. Its value group is denoted by $v_P^{ }F$, and its residue
field by $FP$. When we write $(F|K,P)$ then we mean an extension of
valued fields, that is, $P$ is a place on $L$, and $K$ is endowed with
its restriction (which we will also denote by $P$). This restriction
need not be the identity.

If $F|K$ is a function field and $P$ is a place of $F|K$, then
$K\subseteq FP$. By the \bfind{dimension of $P$} we mean the
transcendence degree $\trdeg FP|K$. Hence, $P$ is called
\bfind{zero-dimensional} if $FP|K$ is algebraic. We will say that
$(F,P)$ has \bfind{rank 1} if $v_P^{ }F$ is archimedean ordered, that
is, embeddable in the ordered additive group of the reals. For the
general definition of the \bfind{rank}, see Section~\ref{sectrank}.

\begin{theorem}                             \label{MTr1E}
Assume that $(F,P)$ has rank 1, and take a subextension $E|K$ of
$F|K$ such that $F|E$ is separable-algebraic. Suppose that $P$
is zero-dimensional. Then in addition to the assertion of
Theorem~\ref{MT1}, ${\cal F}$ can be chosen such that ${\cal F}
|E.{\cal K}$ and ${\cal F}|F.{\cal K}$ are Galois extensions.
\end{theorem}
(Here, ``$F.{\cal K}$'' denotes the field compositum of ${\cal K}$ and
$F$ inside of ${\cal F}$, i.e., the smallest subfield of ${\cal F}$
containing ${\cal K}$ and $F$.)

For rank 1 places of non-zero dimension, we can still prove:
\begin{theorem}                             \label{MTr12}
Assume that $(F,P)$ has rank 1. Then in addition to the assertion of
Theorem~\ref{MT1}, ${\cal F}$ can be obtained from $F.{\cal K}$ by at
most two Galois extensions.
\end{theorem}

Next, we will discuss some important special cases, in particular those
where we can uniformize without extending the function field. For every
place $P$ of $F|K$, we have the following inequality (we will introduce
a more general inequality (\ref{wtdgeq}) later):
\begin{equation}                            \label{Abhie}
\trdeg F|K \>\geq\> \trdeg FP|K \,+\, \dim_\Q \Q\otimes v_P^{ }F\;.
\end{equation}
This is a special case of the Abhyankar inequality. Note that $\dim_\Q
\Q\otimes v_P^{ }F$ is the \bfind{rational rank} of the value group
$v_P^{ }F$, i.e., the maximal number of rationally independent elements
in $v_P^{ }F$. We call $P$ an \bfind{Abhyankar place} if equality holds
in (\ref{Abhie}). An arbitrary place $P$ of a function field $F|K$ is
called \bfind{rational} if $FP=K$.
\begin{theorem}                             \label{MTr1A}
Assume that $P$ is an Abhyankar place of $F|K$ and that $(F,P)$ has rank
1. Then ${\cal F}$ can always be obtained from $F.{\cal K}$ by a single
Galois extension, and the following additional assertions hold:
\sn
a) \ If $FP|K$ is separable, then we can choose ${\cal K}=K$.
\n
b) \ If $P$ is zero-dimensional, then there is a finite Galois extension
${\cal K}'|{\cal K}$ such that we can set ${\cal F}=F.{\cal K}'$, that
is, ${\cal F}$ can be obtained from $F$ by a normal constant extension.
\n
c) \ If $P$ is rational, then we can choose ${\cal F}=F$ and
${\cal K}=K$.
\end{theorem}

Let us note that if $FP$ can be embedded over $K$ in a trivially valued
subfield of $F$, then we may replace $K$ by the image of $FP$. In this
way, this more general case is subsumed under the case of rational
places. For example, if $FP|K$ is a rational function field, then $FP$
can always be embedded in such a way. On the other hand, if $FP|K$ is
separable, then $FP$ can always be embedded in the henselization of
$(F,P)$: choose a separating transcendence basis, embed the rational
function field generated by it in a trivially valued subfield of $F$,
and extend the embedding by Hensel's Lemma. As an algebraic extension of
a trivially valued field is again trivially valued, the same will hold
for the image of $FP$ in $F$. Hence if $FP|K$ is a function field, then
we only have to take a finite extension of $F$ within its henselization
in order to reduce to the case of rational places.

\pars
There is yet another interesting particular case. At first sight, it
seems to be completely opposed to the case of Abhyankar places; but see
Theorem~\ref{MTdens} below and the subsequent remark. A valued field
$(F,P)$ is called \bfind{discretely valued} and $P$ is called
\bfind{discrete} if $v_P^{ }F\simeq\Z$.

\begin{theorem}                                 \label{MTdisc}
Assume that $P$ is a rational discrete place of $F|K$.
Then we can choose ${\cal F}=F$ and ${\cal K}=K$.
\end{theorem}

The assertions of this theorem and of part c) of Theorem~\ref{MTr1A}
can also be formulated as follows:
\begin{corollary}
On a function field $F|K$, all rational rank 1 Abhyankar places
and all rational discrete places are uniformizable.
\end{corollary}

\pars
In [K3], we show that the zero-dimensional rank 1 Abhyankar places, as
well as the zero-dimensional discrete places, lie dense in the Zariski
space of all places of $F|K$, with respect to a ``Zariski patch
topology''. This topology is finer than the Zariski topology (but still
compact); its basic open sets are the sets of the form
\[\{P\mid \mbox{$P$ a place of $F|K$ such that }a_1 P\ne 0,\ldots,
a_k P\ne 0\,;\,b_1 P=0,\ldots,b_\ell P=0\}\]
with $a_1 ,\ldots,a_k,b_1,\ldots, b_\ell\in F\setminus \{0\}$.

\pars
If $K$ is algebraically closed, then every zero-dimensional place of
$F|K$ is rational. So Theorem~\ref{MTr1A} shows (and the first assertion
also follows from Theorem~\ref{MTdisc}):
\begin{corollary}
If $K$ is algebraically closed, then the uniformizable places of $F|K$
lie dense in the Zariski space of~$F|K$, with respect to the Zariski
patch topology. If $K$ is perfect, then the same holds for the places of
$F|K$ which are uniformizable with respect to given elements of $F$
after a finite Galois constant extension of $F$.
\end{corollary}

\parb
Now we turn to places of arbitrary rank.
\begin{theorem}                             \label{MTr>AG}
Assume that $P$ is an Abhyankar place of rank $r>1$ of the function
field $F|K$. Then we can obtain ${\cal F}$ from $F.{\cal K}$ by a
sequence of at most $r-1$ Galois extensions if $P$ is zerodimensional,
or at most $r$ Galois extensions otherwise.
\end{theorem}

In [K6] we will show that $(F|K,P)$ is always weakly uniformizable if
$P$ is an Abhyankar place of $F|K$ for which $FP|K$ is separable. But we
do not obtain that they are uniformizable in general.

For non-Abhyankar places of arbitrary rank, we are not able to prove
that one can obtain ${\cal F}$ by Galois extensions. The obstruction
is, roughly speaking, that we work with extensions in henselizations,
but that taking normal hulls of such extensions may lead to inseparable
residue field extensions. But we will show in [K6] that ${\cal F}$ can
be taken such that it differs from a Galois extension of $F.{\cal K}$
only by an extension in the henselization.

The construction of places given in [K3] yields Abhyankar places or,
if so desired, non-Abhyankar places which are still ``very close to''
Abhyankar places: they lie in the completion of a subfield on which
their restriction is an Abhyankar place. Therefore, it is important to
know that also such places are uniformizable. By ``completion'' we mean
the completion with respect to the uniformity induced by the valuation.
Note that $(F',P)$ lies in the completion of $(F,P)$ if it is an
extension of $(F,P)$ satisfying that for every $a\in F'$ and $\alpha\in
v_P^{ }F'$ there is some $b\in F$ such that $v_P^{ }(a-b)\geq\alpha$.

\begin{theorem}                                  \label{MTdens}
If $(F|K,P)$ satisfies the assumptions of Theorem~\ref{MTr1A} or
Theorem~\ref{MTr>AG}, then the assertions of these theorems carry
over to every function field $(F'|K,P)$ for which $(F',P)$ lies in the
completion of $(F,P)$.
\end{theorem}

\begin{remark}
Theorem~\ref{MTdisc} follows directly from this theorem together with
part c) of Theorem~\ref{MTr1A}. Indeed, if $P$ is a rational discrete
place of $F|K$ and we choose $x\in F$ such that $v_P^{ }x$ is the
smallest positive element in $v_P^{ }F$, then (the restriction of) $P$
is a rational Abhyankar palce of $K(x)|K$ and $F$ lies in the completion
of $(K(x),P)$.
\end{remark}

\parb
In characteristic 0, Theorem~\ref{MT1} is obviously weaker than
Zariski's original result. On the other hand, our proof will yield an
interesting additional assertion. In general, it seems impossible
to obtain it without taking into the bargain a finite extension of the
function field (see the example given in [K5], [K6]).
Let us consider a place $P$ of the function field $E|K$.
We set $\rho=\dim_\Q \Q\otimes v_P^{ }E$ and $\tau= \trdeg EP|K$. We
take elements $x_1,\ldots,x_{\rho}\in E$ such that $v_P^{ }x_1,\ldots,
v_P^{ } x_{\rho}$ are rationally independent elements in $v_P^{ } E$.
Further, we take elements $y_1,\ldots, y_{\tau}\in E$ such that $y_1P,
\ldots, y_{\tau} P$ are algebraically independent over $K$. Then $x_1,
\ldots, x_{\rho},y_1,\ldots,y_{\tau}$ are algebraically independent over
$K$ (cf.~Theorem~\ref{prelBour}) and therefore, $\rho+ \tau\leq \trdeg
E|K$. Every subfield $K(x_1,\ldots,x_{\rho},y_1,\ldots, y_{\tau})$ of
$E$ obtained in this way will be called an \bfind{Abhyankar field of
$(E|K,P)$}. Note that if $P$ is a place of a function field $F|K$ and
$E|K$ is a subextension such that $F|E$ is algebraic, then an Abhyankar
field of $(E|K,P)$ will also be an Abhyankar field of $(F|K,P)$.
\begin{theorem}                             \label{MT4}
Assume the situation as given in Theorems~\ref{MT1} and~\ref{MT1norm}.
Suppose in addition that $(E,P)$ has rank 1, and take any Abhyankar
field $E_0$ of $(E|K,P)$. Then in addition to the assertion of
Theorem~\ref{MT1} and Theorem~\ref{MT1norm}, $\,t_1,\ldots,
t_{\rho+\tau}$ can be chosen algebraic over $E_0\,$.
\end{theorem}
\n
We do not know whether this theorem can be proved for arbitray rank
if one insists in taking normal or Galois extensions. See [K5] and
[K6] for other versions.

%
%
%
%

\parb
To describe a necessary condition for ${\cal F}$ to be equal to $F$, we
need some further definitions. For an arbitrary valued field $(F,P)$ and
a given extension of $P$ from $F$ to its separable-algebraic closure
$F\sep$, the \bfind{absolute inertia field} is defined to be the inertia
field of the normal extension $(F\sep|F,P)$. The decomposition field of
$(F\sep|F,P)$ is the \bfind{henselization of $(F,P)$ in $(F\sep,P)$};
we will denote it by $(F^h,P)$. A valued function field $(F|K,P)$ will
be called \bfind{inertially generated} if it admits a transcendence
basis $T$ such that $(F,P)$ lies in the absolute inertia field of
$(K(T),P)$ (for some extension of $P$ from $K(T)$ to $K(T)\sep$). If it
admits a transcendence basis $T$ such that $(F,P)$ lies in the
henselization of $(K(T),P)$, then we call it \bfind{henselian
generated}.
\begin{theorem}                             \label{MT5}
Assume that $({\cal F}|{\cal K},P)$ is weakly uniformizable (where $P$
is not necessarily the identity on ${\cal K}$). Then $({\cal F}|{\cal
K},P)$ is inertially generated. In particular, ${\cal F}|{\cal K}$ is
separable. If in addition ${\cal F}P={\cal K}$, then $({\cal F}|
{\cal K},P)$ is even henselian generated.
\end{theorem}
\n
{\bf OPEN PROBLEM:} \ Is every inertially
generated valued function field weakly uniformizable?

\bn
\parb
We will deduce Theorems~\ref{MT1} through~\ref{MT4} from two main
theorems which we proved in [K1] (cf.\ also [K2]).

The first theorem is a generalization of the ``Grauert--Remmert
Stability Theorem''. To state it, we introduce a fundamental notion.
Every finite extension $(L|K,P)$ of valued fields satisfies the {\bf
fundamental inequality} (cf.\ [EN], [R], [Z--S] or [K2]):
\begin{equation}                             \label{fundineq}
n\>\geq\>\sum_{i=1}^{\rm g} {\rm e}_i {\rm f}_i
\end{equation}
where $n=[L:K]$ is the degree
of the extension, $P_1,\ldots,P_{\rm g}$ are the distinct extensions of
$P$ from $K$ to $L$, ${\rm e}_i=(v_{P_i}L:v_P^{ }K)$ are the respective
ramification indices and ${\rm f}_i=[LP_i:KP]$ are the respective
inertia degrees. Note that ${\rm g}=1$ if $(K,P)$ is henselian.

A valued field $(K,P)$ is called \bfind{defectless} (or \bfind{stable})
if equality holds in (\ref{fundineq}) for every finite extension $L|K$.
If $\chara KP=0$, then $(K,P)$ is defectless (this is a consequence of
the ``Lemma of Ostrowski'', cf.\ [EN], [R], [K2]).

Now let $(L|K,P)$ be any extension of valued fields (we do not
require that $P$ be the identity on $K$). Assume that $L|K$ has finite
transcendence degree. Then (by Corollary~\ref{fingentb} below):
\begin{equation}                            \label{wtdgeq}
\trdeg L|K \>\geq\> \trdeg LP|KP \,+\, \dim_\Q \Q\otimes
(v_P^{ }L/v_P^{ }K)\;.
\end{equation}
%
%
Note that (\ref{Abhie}) is the special version of (\ref{wtdgeq}) for the
case of a valued function field with trivially valued base field. We
will say that $(L|K,P)$ is \bfind{without transcendence defect} if
equality holds in (\ref{wtdgeq}). Now we are able to state the
\bfind{Stability Theorem}, which deals with an arbitrary valued
function field $(F|K,P)$ (where $P$ need not be the identity on $K$).
\begin{theorem}                           \label{ai}
Let $(F|K,P)$ be a valued function field without transcendence defect.
If $(K,P)$ is a defectless field, then also $(F,P)$ is a defectless
field.
\end{theorem}

The second theorem is a structure theorem for immediate function fields.
An extension $(L|K,P)$ is called \bfind{immediate} if the canonical
embedding of $v_P^{ }K$ in $v_P^{ }L$ and the canonical embedding of
$KP$ in $LP$ are surjective (we then write $v_P^{ }K=v_P^{ }L$ and
$KP=LP$). In this paper, we will only need a special case of the
theorem. We will state and employ the full theorem in [K6].
\begin{theorem}                \label{stt3}
Let $K$ be a separable-algebraically closed field and $(F|K,P)$ an
immediate function field of trans\-cendence degree 1. If $F|K$ is
separable, then
\begin{equation}
\mbox{there is $x\in F$ such that }\; (F^h,P)\,=\,(K(x)^h,P)\;,
\end{equation}
that is, $(F|K,P)$ is henselian generated.
\end{theorem}
For valued fields of residue characteristic 0, the assertion is a direct
consequence of the fact that every such field is defectless (in fact,
every $x\in F\setminus K$ will then do the job). In contrast to this,
the case of positive residue characteristic requires a much deeper
structure theory of immediate algebraic extensions of henselian fields,
in order to find suitable elements $x$.

\parm
In Chapter~\ref{sectprel}, we introduce some further valuation
theoretical tools, including a part of Kaplansky's theory of
immediate extensions, which will also play a crucial role in our proofs.
In Chapter~\ref{sect3}, we give a criterion for valued function
fields to be inertially generated, and prove Theorem~\ref{MT5}. In
Chapter~\ref{sectrup}, we prove  Theorem~\ref{translu}.
This transitivity result allows us to build up our function fields
by various sorts of algebraic and transcendental extensions which all
can be shown separately (in Chapter~\ref{sectufe}) to be uniformizable.
Finally, in Chapters~\ref{sectGnu}, \ref{sectprMT1} and~\ref{sectprMT}
we put everything together to prove our main theorems.

\bn
\parb
I would like to thank Peter Roquette, Mark Spivakovsky,
Bernard Teissier and Frans Oort for support and inspiring
conversations, and the staff of the Fields Institute for their
hospitality. Very special thanks to Hans Schoutens for many
pleasant and encouraging discussions.

%
%
\section{Valuation theoretical preliminaries}         \label{sectprel}
For basic facts from valuation theory, see [EN], [R], [W], [Z--S], [K2].

We will denote the algebraic closure of a field $K$ by $\tilde{K}$.
Whenever we have a place $P$ on $K$, we will automatically fix
an extension of $P$ to the algebraic closure $\tilde{K}$ of $K$. It does
not play a role which extension we choose, except if $P$ is also given
on an extension field $L$ of $K$; in this case, we choose the extension
to $\tilde{K}$ to be the restriction of the extension to $\tilde{L}$. We
say that $P$ is \bfind{trivial} on $K$ if it is an isomorphism of $K$,
which is equivalent to $v_P^{ }K=\{0\}$. If $P$ is given on some
extension field $L$ of $K$ and is trivial on $K$, then there is some
place $P'$ of $L$ which is equivalent to $P$ (i.e., they have the same
valuation ring on $L$) and whose restriction to $K$ is the identity.

A valued field is \bfind{henselian} if it satisfies Hensel's Lemma; see
[R], [W], [K2]. Originally, Hensel's Lemma was proved for complete
discretely valued fields. But it also holds for complete valued fields
of rank~1; i.e., such fields are henselian (cf.\ [W], Theorem 32.11, or
[K2]). The henselization $(K^h,P)$ of a valued field $(K,P)$ (in
$(K\sep,P)\,$) is the minimal henselian extension of $(K,P)$, in the
following sense: if $(L,P')$ is a henselian extension field of $(K,P)$,
then there is a unique embedding of $(K^h,P)$ in $(L,P')$. This is the
\bfind{universal property of the henselization}. We note that every
algebraic extension of a henselian field is again henselian. In
particular, since the absolute inertia field of an arbitrary valued
field contains its henselization, it is henselian.

The following lemma is proved in [K2] (and partially also in [EN]):
\begin{lemma}                               \label{hdl}
A valued field $(K,P)$ is defectless if and only if its henselization
$(K^h,P)$ is.
\end{lemma}

\begin{corollary}                           \label{corhdie}
If $(K,P)$ is defectless, then $(K^h,P)$ does not admit proper immediate
algebraic extensions.
\end{corollary}
\begin{proof}
If $(K,P)$ is defectless, then so is $(K^h,P)$, by the foregoing lemma.
Suppose that $(L|K^h,P)$ is a finite immediate algebraic extension.
Hence, $(v_P^{ }L:v_P^{ }K^h)=1=[LP:K^hP]$. Since $(K^h,P)$ is a
henselian field, there is a unique extension of $v_P^{ }$ from $K^h$ to
$L$. Since $(K^h,P)$ is defectless, we have that $[L:K^h]=(v_P^{ }L:
v_P^{ }K^h)[LP:K^hP]=1$, showing that $L=K^h$. As every proper immediate
extension would contain a proper finite immediate extension, it follows
that $(K^h,P)$ does not admit any proper immediate algebraic extension.
\end{proof}

We also note (see [K2] for the easy proof):
\begin{lemma}                               \label{KsacP}
If $K$ is an arbitrary field and $P$ is a place on $K\sep$, then
$v_P^{ } K\sep$ is the divisible hull of $v_P^{ }K$, and $(KP)\sep
\subseteq K\sep P$. If in addition $P$ is non-trivial on $K$, then
$K\sep P$ is the algebraic closure of $KP$.
\end{lemma}


\begin{lemma}                               \label{ifcont}
Let $P$ be a place of $F|K$ and suppose that $E$ is a subfield of $F$
on which $P$ is trivial. Let $(F^i,P)$ denote the absolute inertia field
of $(F,P)$. Then $E\sep\subset F^i$. Further, if $FP|EP$ is algebraic,
then $(F.E\sep)P$ is the separable-algebraic closure of $FP$.
\end{lemma}
\begin{proof}
By assumption, $P$ induces an embedding of $E$ in $FP$. Further,
we know by ramification theory ([EN], [K2]) that $F^iP$ is
separable-algebraically closed. Thus, $(EP)\sep\subset F^iP$. Using
Hensel's Lemma, one shows that the inverse of the isomorphism $P|_E$ can
be extended from $EP$ to an embedding of $(EP)\sep$ in $F^i$. Its image
is separable-algebraically closed and contains $E$. Hence, $E\sep\subset
F^i$. Further, $(F.E\sep)P$ contains $E\sep P$, which by
Lemma~\ref{KsacP} contains $(EP)\sep$. As $F.E\sep|F$ is algebraic,
so is $(F.E\sep)P|FP$. Therefore, if $FP|EP$ is algebraic, then
$(F.E\sep)P$ is algebraic over $(EP)\sep$ and hence
separable-algebraically closed. Since $(F.E\sep)P\subset
F^iP=(FP)\sep$, it follows that $(F.E\sep)P=(FP)\sep$.
\end{proof}

\begin{lemma}
Take any extension $L|K$ of valued fields and assume that $P$ is a place
of $L$ which is trivial on $K$. If $LP|K$ is separable, then also $L|K$
is separable.
\end{lemma}
\begin{proof}                               \label{rfsepffsep}
If $LP|K$ is separable and $P'$ is equivalent to $P$, then also $LP'|K$
is separable; thus, we can assume that the restriction of $P$ to $K$ is
the identity. Take a finite purely inseparable extension $K'|K$; we have
to show that it is linearly disjoint from $L|K$. As $P$ is the identity
on $K$ and $K'|K$ is purely inseparable, $P$ is also the identity on
$K'$. Hence, $K'\subset (L.K')P$. It follows that
\[[K':K]\>\geq\> [L.K':L]\>\geq\> [(L.K')P:LP]\>\geq\>[LP.K':LP]
\>=\> [K':K]\;,\]
where the last equality holds since $LP|K$ is separable by assumption.
Hence, equality must hold everywhere, showing that $[K':K]=[L.K':L]$,
i.e., $K'|K$ is linearly disjoint from $L|K$.
\end{proof}
\n
A generalization of this lemma to the case of $P$ not being trivial on
$K$ is stated in [K1].


%
%
\subsection{Valuation independence}
For the easy proof of the following theorem, see [B], Chapter VI,
\S10.3, Theorem~1, or [K2].
\begin{theorem}                                \label{prelBour}
Let $(L|K,P)$ be an extension of valued fields. Take elements $x_i,y_j
\in L$, $i\in I$, $j\in J$, such that the values $v_P^{ }x_i\,$, $i\in
I$, are rationally independent over $v_P^{ }K$, and the residues $y_jP$,
$j\in J$, are algebraically independent over $KP$. Then the elements
$x_i,y_j$, $i\in I$, $j\in J$, are algebraically independent over $K$.

Moreover, if we write
\[f\>=\> \displaystyle\sum_{k}^{} c_{k}\,
\prod_{i\in I}^{} x_i^{\mu_{k,i}} \prod_{j\in J}^{} y_j^{\nu_{k,j}}\in
K[x_i,y_j\mid i\in I,j\in J]\]
in such a way that for every $k\ne\ell$
there is some $i$ s.t.\ $\mu_{k,i}\ne\mu_{\ell,i}$ or some $j$ s.t.\
$\nu_{k,j}\ne\nu_{\ell,j}\,$, then
\[v_P^{ }f\>=\>\min_k\, v_P^{ }\,c_k \prod_{i\in I}^{}
x_i^{\mu_{k,i}}\prod_{j\in J}^{} y_j^{\nu_{k,j}}\>=\>
\min_k\left(v_P^{ }\,c_k\,+\,\sum_{i\in I}^{} \mu_{k,i} v_P^{ }
x_i\right)\;.\]
That is, the value of the polynomial $f$ is equal to the least of the
values of its monomials. In particular, this implies:
\begin{eqnarray*}
v_P^{ }K(x_i,y_j\mid i\in I,j\in J) & = & v_P^{ }K\oplus\bigoplus_{i\in I}
\Z v_P^{ }x_i\\
K(x_i,y_j\mid i\in I,j\in J)P & = & KP\,(y_jP\mid j\in J)\;.
\end{eqnarray*}
It also implies that the valuation $v_P^{ }$ on $K(x_i,y_j\mid i\in
I,j\in J)$ is uniquely determined by its restriction to $K$, the values
$v_P^{ }x_i$ and the residues $y_jP$.
\end{theorem}

\begin{corollary}                              \label{fingentb}
Let $(L|K,P)$ be an extension of valued fields of finite transcendence
degree. Then (\ref{wtdgeq}) holds. If in addition $L|K$ is a function
field, and if equality holds in (\ref{wtdgeq}), i.e., $(L|K,P)$ is a
valued function field without transcendence defect, then the extensions
$v_P^{ }L| v_P^{ }K$ and $LP|KP$ are finitely generated. In particular,
if $P$ is trivial on $K$, then $v_P^{ }L$ is a product of finitely many
copies of $\Z$, and $LP$ is again a function field over $K$.
\end{corollary}
\begin{proof}
Choose elements $x_1,\ldots,x_{\rho},y_1,\ldots,y_{\tau}\in L$ such
that the values $v_P^{ }x_1,\ldots,v_P^{ }x_{\rho}$ are rationally
independent over $v_P^{ }K$ and the residues $y_1P,\ldots,y_{\tau} P$
are algebraically independent over $KP$. Then by the foregoing lemma,
$\rho+\tau\leq\trdeg L|K$. This proves that $\trdeg LP|KP$ and the
rational rank of $v_P^{ }L/v_P^{ }K$ are finite. Therefore, we may
choose the elements $x_i,y_j$ such that $\tau=\trdeg LP|KP$ and
$\rho=\dim_{\Q} \Q\otimes (v_P^{ }L/v_P^{ }K)$ to obtain inequality
(\ref{wtdgeq}).

Assume that this is an equality. This means that for $L_0:=K(x_1,\ldots,
x_{\rho},y_1,\ldots,y_{\tau})$, the extension $L|L_0$ is algebraic.
Since $L|K$ is finitely generated, it follows that $L|L_0$ is finite. By
the fundamental inequality, this yields that $v_P^{ }L|v_P^{ }L_0$ and
$LP| L_0P$ are finite extensions. Since already $v_P^{ } L_0|v_P^{ } K$
and $L_0P|KP$ are finitely generated by the foregoing theorem, it
follows that also $v_P^{ }L|v_P^{ }K$ and $LP|KP$ are finitely
generated.
\end{proof}

\begin{lemma}                               \label{PQAbh}
Take a function field $F|K$ and a place $P$ of $F|K$. If $P=Q\ovl{Q}$,
then
\[\dim_\Q \Q\otimes v_P^{ }F\>=\> \dim_\Q \Q\otimes v_Q^{ }F\,+\,
\dim_\Q \Q\otimes v_{\ovl{Q}} FQ\;.\]
Further, $P$ is an Abhyankar place if and only if $Q$ and $\ovl{Q}$ are.
\end{lemma}
\begin{proof}
The first assertion is well-known. The second assertion follows
from the first, keeping in mind that $FP=(FQ)\ovl{Q}$. We leave
the straightforward proofs to the reader.
\end{proof}

\subsection{The rank}                       \label{sectrank}
The \bfind{rank of an ordered abelian group} $\Gamma$ is the order
type of the chain of its proper convex subgroups. By a theorem of
H\"older, $\Gamma$ is embeddable in the ordered additive group of the
reals if and only if its rank is 1, that is, the only proper
convex subgroup is $\{0\}$. This in turn holds if and only if
$\Gamma$ is \bfind{archimedean ordered}, i.e., for every two positive
elements $\alpha,\beta\in\Gamma$ there is some $n\in\N$ such that
$n\alpha\geq\beta$.

The \bfind{rank of} a valued field $(F,P)$ is defined to be the rank of
its value group $v_P^{ }F$. If this is finite, say $n$, then $P$ is the
composition of $n$ places: $P=P_1P_2\ldots P_n\,$, where all
$P_i$ have value groups of rank 1 (cf.\ [Z--S], [K2]).

If $\Gamma_1$ is a subgroup of $\Gamma$, then its divisible hull $\Q
\otimes\Gamma_1$ lies in the convex hull of $\Gamma_1$ in $\Q\otimes
\Gamma$. Hence if $\Gamma_1$ is a proper convex subgroup of $\Gamma$,
then $\Q\otimes\Gamma_1$ is a proper convex subgroup of $\Q\otimes
\Gamma$ and thus, $\dim_\Q \Q\otimes\Gamma_1<\dim_\Q \Q\otimes\Gamma$.
It follows that if $\{0\}=\Gamma_0\subsetuneq\Gamma_1\subsetuneq\ldots
\subsetuneq \Gamma_n= \Gamma$ is a chain of convex subgroups of
$\Gamma$, then $\dim_\Q \Q\otimes\Gamma\geq n$. In view of
(\ref{Abhie}), this proves that the rank of a place $P$ of a function
field $F|K$ cannot exceed $\trdeg F|K$ and thus is finite.

If $(F'|F,P)$ is an algebraic extension of valued fields, then
$v_P^{ }F'/v_P^{ }F$ is a torsion group and $F'P|FP$ is algebraic
(this is a consequence of Theorem~\ref{prelBour}). In particular,
the rank of $v_PF'$ is equal to that of $v_PF$. Hence, passing to
an algebraic extension does not change the rank of a valued field. As
the rank of an ordered abelian group does not increase by passing to a
subgroup, the rank of a valued field does not increase by passing to a
subfield.

%
%
\subsection{Kaplansky approximation}
The material of this section is based on work by Ostrowski and
Kaplansky [KA] (cf.\ also [K2]). The proof of the first lemma is an
easy exercise.
\begin{lemma}                               \label{nomax}
The extension $(L|K,P)$ is immediate if and only if for every $z\in L$,
the set $\{v_P^{ }(z-a)\mid a\in K\}$ has no maximal element.
\end{lemma}
\begin{lemma}                               \label{close}
Let $(K(z)|K,P)$ be an immediate transcendental extension. Assume that
$(K,P)$ is a separable-algebraically closed field or that $(K(z),P)$
lies in the completion of $(K,P)$. Take any polynomial $f\in K[X]$. Then
the value $v_P^{ }f(a)$ is fixed for all $a\in K$ sufficiently close to
$z$. That is,
\begin{equation}                            \label{trat}
\begin{array}{c}
\forall f\in K[X]\;\exists\alpha\in v_P^{ }K\;\exists\beta\in
\{v_P^{ }(z-b)\mid b\in K\}\; \forall a\in K:\\[.3cm]
v_P^{ }(z-a) \geq\beta\,\Rightarrow\, v_P^{ }f(a)=\alpha\;.
\end{array}\;\;\;
\end{equation}
\end{lemma}
Kaplansky proves that if (\ref{trat}) does not hold, then there is a
proper immediate algebraic extension of $(K,P)$. If $(K(z),P)$ does not
lie in the completion of $(K,P)$, then this can be transformed into a
proper immediate separable-algebraic extension ([K1], [K2]; the proof
uses a variant of the Theorem on the Continuity of Roots). But such
an extension cannot exist if we assume that $K$ be
separable-algebraically closed.
If on the other hand $(K(z),P)$ lies in the completion of
$(K,P)$, then one can show that if $f$ does not satisfy (\ref{trat}),
then $v_P^{ }f(z)=\infty$. But this means that $f(z)=0$, contradicting
the assumption that $K(z)|K$ is transcendental.

\pars
For a polynomial $f$ in one variable over a field of arbitrary
characteristic, we denote by $f\iT $ its $i$-th formal derivative
(cf.\ [KA], [K2]). These polynomials are defined such that the following
Taylor expansion holds in arbitrary characteristic:
\[f(z)\>=\>f(a)+\sum_{i=1}^{\deg f} f\iT (a)(z-a)^i\;.\]

\begin{lemma}                               \label{lvpol}
Assume that (\ref{trat}) holds, and take any polynomial $f\in K[z]$.
Then there are $a,b\in K$ such that for $\tilde{z}:= \frac{z-a}{b}\in
K[z]$, we have that $v_P^{ }\tilde{z}=0$ and the values of the non-zero
among the elements $f\iT (a)b^i\in K$ are all distinct. With such $a$
and $b$,
\begin{equation}                            \label{vpol}
v_P^{ } f(z)\>=\>v_P^{ }\sum_{i=0}^{\deg f} f\iT (a)(z-a)^i\>=\>v_P^{ }
\sum_{i=0}^{\deg f} f\iT (a)b^i \tilde{z}^i\>=\>
\min_i\, v_P^{ } f\iT (a)b^i\;.
\end{equation}
If finitely many polynomials in $K[z]$ are given, then $a,b$ can be
chosen such that (\ref{vpol}) holds simultaneously for all of them.
\end{lemma}
\begin{proof}
Take finitely many polynomials $f_1,\ldots,f_n\in K[z]$. From
Lemma~\ref{close} we know that for all $a\in K$ close enough to $z$,
the values $v_P^{ } f_j\iT (a)$ of the non-zero among the polynomials
$f_j\iT$, $i,j\in\N$, are fixed. Since by Lemma~\ref{nomax} the set
$\{v_P^{ }(z-a)\mid a\in K\}$ has no maximal element, we can then take
$a$ so close to $z$ that for every fixed $j$, the values of all non-zero
elements $f_j\iT (a)(z-a)^i$, $i\in\N$, are distinct. Having picked such
an element $a\in K$, we choose an element $b\in K$ such that $v_P^{ }b=
v_P^{ }(z-a)$. Then (\ref{vpol}) holds by the ultrametric triangle law.
\end{proof}


%
%
\subsection{Transcendence bases of separable valued function fields}
We will denote the algebraic closure of $K$ by $\tilde{K}$.
We assume that $K$ is a field and that $P$ is a place on the rational
function field $\tilde{K}(z)$ and infer the following two lemmata from
[KH--K]:
\begin{lemma}
The following assertions are equivalent:\sn
a) \ $(\tilde{K}(z)|\tilde{K},v_P^{ })$ is immediate,\n
b) \ $v_P^{ }K(z)/v_P^{ }K$ is a torsion group and $K(z)P|KP$ is
algebraic,\n
c) \ $\{v_P^{ }(z-c)\mid c\in\tilde{K}\}$ has no maximal element.
\end{lemma}
\begin{lemma}
Assume that $\{v_P^{ }(z-c)\mid c\in\tilde{K}\}$ has a maximal element
and that $c_0\in \tilde{K}$ is an element of minimal degree over $K$
such that $v_P^{ }(z-c_0)$ is such a maximal element. Take $f$ to be the
minimal polynomial of $c_0$ over $K$.
\sn
1) \ If $v_P^{ }K(z)/v_P^{ }K$ is not a torsion group, then $v_P^{ }
f(z)$ is not a torsion element modulo $v_P^{ }K$.\n
2) \ If $K(z)P|KP$ is transcendental, then there is some $e\in\N$ and
some $d\in K$ such that $(df(z)^e)P$ is transcendental over $KP$.
\end{lemma}

From these we deduce:
\begin{lemma}
Assume that $(\tilde{K}(z)|\tilde{K},v_P^{ })$ is not immediate. Then
there is some $h\in K[X]$ such that $v_P^{ }h(z)$ is non-torsion
over $v_P^{ }K$ or that $h(z)P$ is transcendental over $KP$, and such
that $K(z)|K(h(z))$ is separable.
\end{lemma}
\begin{proof}
If $v_P^{ }z$ is not torsion modulo $v_P^{ }K$ or $zP$ is transcendental
over $KP$, then we set $h(X):=X$.

Otherwise, we set $g(X):=f(X)$ if case 1) of the foregoing lemma holds,
and $g(X):=df(X)^e$ if case 2) holds (by Theorem~\ref{prelBour}, only one
of the two cases can hold at a time). If the polynomial $g(X)-g(z)$ is
separable over $K(g(z))$, then we set $h(X):=g(X)$. Otherwise, we
proceed as follows. Set $n:=\deg g$; this must be divisible by the
characteristic $p$ of $K$. In case 1), $v_P^{ }z$ is torsion modulo
$v_P^{ }K$ by assumption, and it follows that for $h(X):= Xg(X)$, the
value $v_P^{ }h(z)=v_P^{ }z+v_P^{ }g(z)$ is still not torsion modulo
$v_P^{ }K$. As $\deg h(X) =n+1$ is not divisible by $p$, we find that
$h(X)-h(z)$ is separable over $K(h(z))$. In case 2), $v_P^{ }z$ is
torsion modulo $v_P^{ }K$ by (\ref{wtdgeq}), and it follows that there
is some $b\in K$ such that $v_P^{ }bz>0$. Then for $h(X):=g(X)+bX$ we
have that $h(z)P=g(z)P$, and that $h(X)-h(z)$ is separable over
$K(h(z))$.
\end{proof}

Now we are able to prove:
\begin{lemma}                               \label{tb}
Let $(F|K,P)$ be any valued function field (where $P$ is not necessarily
trivial on $K$). Assume that $F|K$ is separable. Then there is a
separating transcendence basis of $F|K$ containing elements $x_1,\ldots,
x_{\rho},y_1,\ldots,y_{\tau}$ such that $v_P^{ }x_1,\ldots, x_{\rho}
v_P^{ }$ is a maximal set of elements in $v_P^{ }F$ rationally
independent modulo $v_P^{ }K$, and $y_1P,\ldots,y_{\tau}P$ form a
transcendence basis of $FP|KP$.
\end{lemma}
\begin{proof}
Since $F|K$ is separable, we can choose a separating transcendence basis
$z_1,\ldots,z_n$ of $F|K$. We set $K_0:=K$ and $K_i:=K(z_1,\ldots,z_i)$.
We proceed by induction on $i$. If the extension $(\tilde{K}_{i-1}(z_i)|
\tilde{K}_{i-1},v_P^{ })$ is not immediate, then we choose $h_i(z_i)$
according to the assertion of the foregoing lemma. Otherwise, we set
$h_i(X):=X$. Since every extension $K_{i-1}(z_i)|K_{i-1}(h_i(z_i))$ is
separable, we obtain a separating transcendence basis
$h_1(z_1),\ldots,h_n(z_n)$ of $F|K$.

We set $\rho:=\dim_{\Q}\Q\otimes(v_P^{ }F/v_P^{ }K)$ and $\tau:=\trdeg
FP|KP$. As
\[\rho\,=\,\sum_{i=0}^{n-1}\dim_{\Q}\Q\otimes(v_P^{ }K_{i+1}/
v_P^{ }K_i)\mbox{ \ \ \ and \ \ \ }\tau\,=\, \sum_{i=0}^{n-1} \trdeg
K_{i+1}P|K_iP\;,\]
and in view of the fact that
\[\dim_{\Q}\Q\otimes(v_P^{ }K_{i+1}/v_P^{ }K_i) \,+\,\trdeg
K_{i+1}P|K_iP\>\leq\> \trdeg K_{i+1}|K_i\>=\>1\;,\]
we find that for precisely $\rho$ many values of $i$, $v_P^{ }h_i(z_i)$
will be rationally independent modulo $v_P^{ } K_{i-1}\,$. Collecting
all of these $h_i(z_i)$ and calling them $x_1, \ldots, x_{\rho}$ we thus
obtain that $v_P^{ }x_1 ,\ldots, v_P^{ }x_{\rho}$ is a maximal set of
elements in $v_P^{ }F$ rationally independent modulo $v_P^{ }K$.
Similarly, we find that for precisely $\tau$ many values of $\,i$,
the residues $h_i(z_i)P$ will be transcendental over $K_{i-1}P$.
Collecting all of these $h_i(z_i)$ and calling them
$y_1,\ldots,y_{\tau}$ we thus obtain that $y_1P,\ldots,y_{\tau}P$
form a transcendence basis of $FP|KP$.
\end{proof}

\begin{remark}
We do not know whether in addition to the assertion of the lemma, the
$y_i$ can be chosen such that $y_1P,\ldots,y_{\tau}P$ form a separating
transcendence basis of $FP|KP$.
\end{remark}

%
%
\section{Inertially and henselian generated function
fields}                                              \label{sect3}
\begin{theorem}                             \label{hrwtd}
Assume that $F|K$ is a function field and $P$ an Abhyankar place of
$F|K$ such that $FP|K$ is a separable extension. Then $(F|K,P)$ is
inertially generated. If in addition $FP=K$ or $FP|K$ is a rational
function field, then $(F|K,P)$ is henselian generated. If $v_P^{ }F= \Z
v_P^{ }x_1\oplus \ldots\oplus\Z v_P^{ }x_\rho$ and $y_1P,\ldots,y_\tau
P$ is a separating transcendence basis of $FP|K$, then $T=\{x_1,\ldots,
x_\rho,y_1,\ldots, y_\tau\}$ is a generating transcendence basis, that
is, $(F,P)$ lies in the absolute inertia field of $(K(T),P)$, and if
$FP=K(y_1P,\ldots,y_\tau P)$, then $F\subseteq K(T)^h$.
\end{theorem}
\begin{proof}
By Corollary~\ref{fingentb}, value group and residue field of $(F,P)$
are finitely generated. We choose $x_1,\ldots,x_\rho\in F$ such that
$v_P^{ }F=\Z v_P^{ }x_1\oplus\ldots\oplus\Z v_P^{ }x_\rho$. Since $FP|K$
is a finitely generated separable extension, it is separably generated.
Therefore, we can choose $y_1,\ldots,y_\tau \in F$ such that $FP|K(y_1P,
\ldots,y_\tau P)$ is separable-algebraic ($\tau=\trdeg FP|K$). Now we
can choose some $a\in FP$ such that $FP=K(y_1P,\ldots, y_\tau P,a)$.
Since $a$ is separable-algebraic over $K(y_1P,\ldots, y_\tau P)$, by
Hensel's Lemma there exists an element $\eta$ in the henselization of
$(F,P)$ such that $\eta P=a$ and that the reduction of the minimal
polynomial of $\eta$ over $F_0:=K(x_1,\ldots,x_\rho,y_1, \ldots,
y_\tau)$ is the minimal polynomial of $a$ over $KP(y_1P,\ldots, y_\tau
P)$. Then $\eta$ lies in the absolute inertia field of $F_0\,$. Now the
field $F_0(\eta)$ has the same value group and residue field as $F$, and
it is contained in the henselization $F^h$ of $F\,$. As henselizations
are immediate extensions and the henselization $F_0(\eta)^h$ of
$F_0(\eta)$ can be chosen inside of $F^h$, we obtain an immediate
algebraic extension $(F^h| F_0(\eta)^h,P)$. On the other hand, we
observe that $(K,P)$ is a defectless field since $P$ is trivial on $K$.
By construction, $(F_0|K,P)$ is without transcendence defect, and the
same is true for $(F_0(\eta)|K,P)$ since this property is preserved by
algebraic extensions. Hence we know from Theorem~\ref{ai} that
$(F_0(\eta),P)$ is a defectless field. Now Corollary~\ref{corhdie} shows
that the extension $F^h|F_0(\eta)^h$ must be trivial. Therefore, $F$ is
contained in $F_0(\eta)^h$, which in turn is a subfield of the absolute
inertia field of $F_0$. This shows that $(F|K,P)$ is inertially
generated.

If $FP=K$, then we do not need the elements $y_j$ and $a$. If $FP|K$ is
a rational function field, then we can choose $y_1,\ldots,y_\tau\in F$
such that $FP=K(y_1P,\ldots,y_\tau P)$, and we do not need $a$. In both
cases, we find that $F^h=F_0^h$, which yields that $(F|K,P)$ is
henselian generated.
\end{proof}

\bn
$\bullet$ \ {\bf Proof of Theorem~\ref{MT5}}
\sn
Assume that $({\cal F}|{\cal K},P)$ is weakly uniformizable (where $P$
is not necessarily trivial on ${\cal K}$). Denote by $(L,P)$ the
absolute inertia field of $({\cal K}(t_1,\ldots,t_s),P)$.

First, $\det J_{fP} (\eta_1P,\ldots,\eta_nP)\ne 0$ and the fact that
the $f_iP$ are polynomials over ${\cal K}(t_1,\ldots,t_s)P$ imply that
$\eta_1P,\ldots, \eta_nP$ are separable algebraic over ${\cal K}(t_1,
\ldots,t_s)P$ (cf.\ [L], Chapter X, \S7, Proposition 8). On the other
hand, $LP$ is the separable-algebraic closure of ${\cal K}(t_1,\ldots,
t_s)P$. Therefore, there are elements $\eta'_1, \ldots, \eta'_n$ in $L$
such that $\eta'_iP=\eta_iP$. Since $(L,P)$ is henselian, the
multidimensional Hensel's Lemma (cf.\ [K2], [K7]) now shows the
existence of a common root $(\eta''_1,\ldots,\eta''_n)\in L^n$ of the
$f_i$ such that $\eta''_iP=\eta'_iP=\eta_iP$. But by the uniqueness
assertion of the multidimensional Hensel's Lemma (which also holds in
the algebraic closure $\tilde{L}$ of ${\cal K}(t_1,\ldots,t_s)$), we
find that $(\eta''_1,\ldots,\eta''_n)=(\eta_1,\ldots,\eta_n)$. Hence,
the $\eta_i$ are elements of $L$, which proves that $({\cal F}|{\cal K},
P)$ is inertially generated.

\pars
If we have in addition that $P$ is a rational place, then
$\eta_1P,\ldots,\eta_nP\in {\cal K}$. In this case, we can choose
$\eta'_1,\ldots,\eta'_n$ and $\eta''_1,\ldots,\eta''_n$ already in the
henselization of $({\cal K} (t_1,\ldots,t_s),P)$, which implies that
also $\eta_1,\ldots,\eta_n$ lie in this henselization.          \QED

%
%
\section{Basic properties of relative uniformization}   \label{sectrup}
Recall the definition for ``$({\cal F}|{\cal K},P)$ is uniformizable''
given preceding
to Theorem~\ref{MT1}. If $P$ is not trivial on ${\cal K}$, one may think
of this as ``relative uniformization''. For the case of $P$ a place of
$F|K$, relative uniformization will help us to prove Theorem~\ref{MT1},
as we will build up ${\cal F}$ by a tower of uniformizable finitely
generated extensions of valued fields, starting from ${\cal K}$. In the
sections below, we will consider the different types of extensions
involved in the build-up. Beforehand, we need some easy observations.
First, we observe going-up and going-down of uniformizability through
constant extensions of the function field.
\begin{lemma}                               \label{going}
Let $(L|K,P)$ be an extension of valued fields and $F|K$ a function
field such that $F\subset L$. Take $\zeta_1,\ldots,\zeta_m\in
{\cal O}_F\,$.\sn
a) \ Suppose that $(F|K,P)$ is uniformizable with respect to $\zeta_1,
\ldots,\zeta_m\,$. Take an arbitrary subextension $L'|K$ of $L|K$ such
that $\trdeg F.L'|L'=\trdeg F|K$. Then with the same $t_i,\eta_i,f_i$ as
for $(F|K,P)$, also $(F.L'|L',P)$ is uniformizable with respect to
$\zeta_1, \ldots,\zeta_m\,$.
\sn
b) \ Suppose that $L'|K$ is a subextension of $L|K$ such that $(F.L'|L',
P)$ is uniformizable with respect to $\zeta_1, \ldots,\zeta_m\,$. Then
there is a finitely generated subextension $L_0|K$ of $L'|K$ such
that with the same $t_i,\eta_i,f_i$ as for $(F.L'|L',P)$, also
$(F.L_0|L_0,P)$ is uniformizable with respect to $\zeta_1,\ldots,
\zeta_m\,$.
\sn
c) \ Suppose that $(F.K\sep|K\sep,P)$ is uniformizable with respect to
$\zeta_1, \ldots,\zeta_m\,$. Then there is a finite Galois extension
$K'|K$ such that with the same $t_i,\eta_i,f_i$ as for $(F.K\sep|K\sep
,P)$, also $(F.K'|K',P)$ is uniformizable with respect to
$\zeta_1,\ldots,\zeta_m\,$.
\end{lemma}
Assertion a) follows directly from the definition since the condition
``$\trdeg F.L'|L'=\trdeg F|K$'' guarantees that $\{t_1,\ldots,t_s\}$
remains a transcendence basis of $F.L'|L'$. Note that in the case of
$\trdeg F.L'|L'<\trdeg F|K$ it remains valid if there is a transcendence
basis $T\subset {\cal O}_{K(t_1,\ldots,t_s)}$ of $F.L'|L'$ such that
every $t_i$ is contained in ${\cal O}_K [T]$. For the proof of assertion
b), we just have to do the following. We collect the finitely many
coefficients $c_1,\ldots, c_\ell\in L'$ of all polynomials $f_i\in {\cal
O}_{L'}[t_1,\ldots,t_s, X_1,\ldots,X_n]$. Further, we choose a finitely
generated subextension $L'_0|K$ of $L'|K$ so large that $F.L'_0=L'_0
(t_1,\ldots,t_s,\eta_1, \ldots,\eta_n)$ holds. Then we set $L_0:=
L'_0(c_1,\ldots,c_\ell)$. Part c) is a direct consequence of a)
and b).

\parm
Now we turn to the \bfind{transitivity of relative uniformization}.
Take a finitely generated extension $F|K$ and a finitely generated
subextension $F_0|K$. Further, take elements $\zeta_1,\ldots, \zeta_m\in
{\cal O}_F$ and assume that $(F|F_0,P)$ is uniformizable with respect to
$\zeta_1,\ldots,\zeta_m\,$. So there are elements $t_1, \ldots,
t_{\tilde{s}}\in {\cal O}_F\,$, algebraically independent over $F_0\;$
($\tilde{s}\geq 0$), elements
$\tilde{\eta}_1,\ldots,\tilde{\eta}_{\tilde{n}}\in {\cal O}_F\,$, with
$\zeta_1,\ldots, \zeta_m$ among them, and polynomials $\tilde{f}_i(X_1,
\ldots, X_{\tilde{n}}) \in {\cal O}_{F_0} [t_1,\ldots,t_{\tilde{s}},X_1,
\ldots, X_{\tilde{n}}]$, $1\leq i \leq \tilde{n}$, such that
\pars
$\bullet$ \ $F=F_0(t_1,\ldots,t_{\tilde{s}},\tilde{\eta}_1,\ldots,
\tilde{\eta}_{\tilde{n}})$,\par
$\bullet$ \ for $i<j$, $X_j$ does not occur in $\tilde{f}_i\,$,\par
$\bullet$ \ $\tilde{f}_i(\tilde{\eta}_1,\ldots,\tilde{\eta}_{\tilde{n}})
=0$ for $1\leq i\leq \tilde{n}$, and\par
$\bullet$ \ $(\det J_{(\tilde{f}_1,\ldots, \tilde{f}_{\tilde{n}})}
(\tilde{\eta}_1,\ldots,\tilde{\eta}_{\tilde{n}}))P\ne 0$.
\sn
Now we collect the coefficients of all polynomials $\tilde{f}_i\in
{\cal O}_{F_0} [t_1,\ldots,t_{\tilde{s}},X_1,\ldots,X_{\tilde{n}}]$ and
call them the \bfind{uniformization coefficients of $F,\zeta_1,\ldots,
\zeta_m$ in $F_0$}. Then we extend $P$ to $\tilde{F}$ and take any
elements $\zeta'_1,\ldots, \zeta'_{m'} \in {\cal O}_{\tilde{F}_0}$
which include these uniformization
coefficients. Assume that ${\cal F}_0$ is an algebraic extension of
$F_0(\zeta'_1,\ldots,\zeta'_{m'})$ such that $({\cal F}_0|K,P)$ is
uniformizable with respect to $\zeta'_1,\ldots,\zeta'_{m'}\,$. So there
are elements $t_{\tilde{s}+1},\ldots,t_{s}\in {\cal O}_{{\cal F}_0}\,$,
algebraically independent over $K\;$ ($s\geq\tilde{s}$), elements
$\eta_1,\ldots, \eta_{n'} \in {\cal O}_{{\cal F}_0}\,$, with the
elements $\zeta'_1, \ldots, \zeta'_{m'}$ among them, and polynomials
$f_i(X_1,\ldots,X_{n'})\in {\cal O}_K [t_{\tilde{s}+1},\ldots,t_s,X_1,
\ldots, X_{n'}]$, $1\leq i\leq n'$, such that
\pars
$\bullet$ \ ${\cal F}_0= K(t_{\tilde{s}+1},\ldots,t_s,\eta_1,\ldots,
\eta_{n'})$,\par
$\bullet$ \ for $i<j$, $X_j$ does not occur in $f_i\,$,\par
$\bullet$ \ $f_i(\eta_1,\ldots,\eta_{n'}) =0$ for $1\leq i\leq {n'}$,
and\par
$\bullet$ \ $(\det J_{(f_1,\ldots,f_{n'})}(\eta_1,\ldots,\eta_{n'}))P\ne 0$.
\sn
We observe that the elements $t_1,\ldots,t_s$ are algebraically
independent over $K$. We set $n:=n'+\tilde{n}\,$. Trivially, the
polynomials $f_i\,$, $1\leq i\leq {n'}$, can be viewed as polynomials in
${\cal O}_K[t_1,\ldots,t_s,X_1,\ldots,X_n]$. Note that by our choice of
$\zeta'_1,\ldots,\zeta'_{m'}\,$, all ${\cal O}_{F_0}$-coefficients of
the polynomials $\tilde{f}_i\in {\cal O}_{F_0}[t_1,\ldots,t_{\tilde{s}},
X_1, \ldots, X_{\tilde{n}}]$ appear as some $\eta_j\,$, with $j\in
\{1,\ldots,n'\}$; we may assume that all $\eta_j$ are distinct. For
$1\leq i\leq \tilde{n}$, we obtain the polynomial $f_{n'+i}\in
{\cal O}_K [t_1,\ldots, t_s,X_1, \ldots, X_n]$ from the polynomial
$\tilde{f}_i$ as follows:
\sn
i) \ \ for $1\leq j\leq\tilde{n}$, we replace $X_j$ by $X_{n'+j}\,$,\n
ii) \ if an ${\cal O}_{F_0}$-coefficient of $\tilde{f}_i$ is equal to
$\eta_j\,$, $1\leq j\leq n'$, then we replace it by $X_j\,$.
\sn
Accordingly, we set $\eta_{n'+i}:= \tilde{\eta}_i\,$. Then (U1) and (U2)
hold. For $1\leq i\leq \tilde{n}$, we have that
\[\frac{\partial f_{n'+i}}{\partial X_{n'+i}}(\eta_1,\ldots,\eta_n) =
\frac{\partial \tilde{f}_i}{\partial X_i} (\tilde{\eta}_1,\ldots,
\tilde{\eta}_{\tilde{n}})\;.\]
This shows that the diagonal elements of the lower triangular matrix
$J_{(f_1,\ldots,f_n)} (\eta_1,\ldots,\eta_n)$ are precisely the diagonal
elements of the lower triangular matrices $J_{(f_1,\ldots,f_{n'})}
(\eta_1, \ldots,\eta_{n'})$ and $J_{(\tilde{f}_1, \ldots,
\tilde{f}_{\tilde{n}})}(\tilde{\eta}_1,\ldots,\tilde{\eta}_{\tilde{n}})$.
Consequently, $\det J_{(f_1,\ldots,f_n)} (\eta_1,\ldots,\eta_n)P\ne 0$.
That is, for
\begin{eqnarray*}
F.{\cal F}_0 & = & K(t_{\tilde{s}+1},\ldots,t_{s},\eta_1,\ldots,
\eta_{n'}, t_1,\ldots,t_{\tilde{s}}, \tilde{\eta}_1,\ldots,
\tilde{\eta}_{\tilde{n}})
\\
 & = & K(t_1,\ldots,t_s,\eta_1,\ldots,\eta_n)
\end{eqnarray*}
we have that $(F.{\cal F}_0|K,P)$ is uniformizable with respect to
$\zeta_1,\ldots,\zeta_m\,$. We have proved:
\begin{lemma}                               \label{basictrans}
Take a finitely generated extension $F|K$, a finitely generated
subextension $F_0|K$, and $\zeta_1,\ldots,\zeta_m\in {\cal O}_F\,$.
Assume that
\sn
1) \ $(F|F_0,P)$ is uniformizable with respect to $\zeta_1,\ldots,
\zeta_m\,$, with uniformizing transcendence basis $T_1\,$,
and\n
2) \ there is a finite extension ${\cal F}_0|F_0$ such that $({\cal F}_0
|K,P)$ is uniformizable with respect to the uniformization coefficients
of $F,\zeta_1,\ldots,\zeta_m$ in $F_0\,$, with uniformizing
transcendence basis $T_2\,$.
\sn
Then $(F.{\cal F}_0|K,P)$ is uniformizable with respect to
$\zeta_1,\ldots,\zeta_m\,$, with uniformizing transcendence basis
$T_1\cup T_2\,$.

Consequently, if $(F|F_0,P)$ and $(F_0|K,P)$ are uniformizable, then so
is $(F|K,P)$.
\end{lemma}
This lemma proves Theorem~\ref{translu}. It is the basic form of
transitivity, from which we will also derive the transitivity of the
following two properties. Let $(L|K,P)$ be an arbitrary extension of
valued fields, and $E$ any subfield of $L$. Then we will say
that $(L|K,P)$ has \bfind{(relative) Galois-uniformization over $E$} if
$E.K|K$ is a function field and for every choice of elements $\zeta_1,
\ldots, \zeta_m\in {\cal O}_L$ there is a finite Galois extension
${\cal E}|E$ such that $\zeta_1,\ldots,\zeta_m\in {\cal E}.K$ and
$({\cal E}.K|K,P)$ is uniformizable with respect to $\zeta_1,\ldots,
\zeta_m\>$ (observe that ${\cal E}.K|K$ is a function field by the
assumptions on $E$ and ${\cal E}|E$). Note that this property implies
that $L|E.K$ is a separable-algebraic extension since otherwise, there
is some $\zeta\in {\cal O}_L$ which is not contained in ${\cal E}.K$ for
any Galois extension ${\cal E}$ of $E$. Similarly, we will say that
$(L|K,P)$ has \bfind{(relative) normal-uniformization over $E$} if for
every choice of elements $\zeta_1, \ldots, \zeta_m\in {\cal O}_L$ there
is a finite Galois extension ${\cal E}|E$ and a purely inseparable
subextension ${\cal K}|K$ of ${\cal E}.K|K$ such that $\zeta_1,\ldots,
\zeta_m\in {\cal E}.K$ and $({\cal E}.K|{\cal K},P)$ is uniformizable
with respect to $\zeta_1,\ldots, \zeta_m\,$. This implies that $L|E.K$
is algebraic.

\begin{lemma}                               \label{transitive}
Assume that $(M|L,P)$ and $(L|K,P)$ are extensions of valued fields.
Take any subfield $E$ of $M$ such that $E.K|K$ is a function field.
If $(M|L,P)$ has Galois-uniformization over $E$ and $(L|K,P)$ has
Galois-uniformization over some common subfield $E_0$ of $E$ and $L$,
then $(M|K,P)$ has Galois-uniformization over $E$. An analogous
assertion holds for normal-uniformization over $E$.
\end{lemma}
\begin{proof}
We only prove the first assertion; the proof of the second assertion is
similar. Take $\zeta_1,\ldots,\zeta_m\in {\cal O}_M\,$, and choose a
finite Galois extension ${\cal E}'$ of $E$ such that $\zeta_1,\ldots,
\zeta_m\in {\cal E}'.L$ and that $({\cal E}'.L|L,P)$ is uniformizable
with respect to the $\zeta$'s.

By our assumption on $(L|K,P)$, $L$ is algebraic over $E_0.K\,$. Hence
by part a) of Lemma~\ref{going}, there is a finite subextension $L_0|
E_0.K$ of $L|E_0.K$ such that $({\cal E}'.L_0|L_0,P)$ is uniformizable
with respect to the $\zeta$'s. Take $\zeta'_1,\ldots,\zeta'_{m'}\in
{\cal O}_{L}$ to include generators of $L_0$ over $E_0.K$ and the
uniformization coefficients of ${\cal E}'.L_0, \zeta_1,\ldots,\zeta_m$
in $L_0\,$.

By hypothesis, there is a finite Galois extension ${\cal E}_0$ of
$E_0$ such that $\zeta'_1,\ldots,\zeta'_{m'}\in {\cal E}_0.K$ and
$({\cal E}_0.K|K,P)$ is uniformizable with respect to $\zeta'_1,\ldots,
\zeta'_{m'}\,$. Then ${\cal E}:={\cal E}'.{\cal E}_0$ is a finite Galois
extension of $E$. By construction, we have that ${\cal E}_0.K$ is a
finite extension of $L_0$ and that $({\cal E}'.L_0).({\cal E}_0.K)=
{\cal E}.K\,$. Hence by Lemma~\ref{basictrans}, $({\cal E}.K|K,P)$ is
uniformizable with respect to $\zeta_1,\ldots,\zeta_m\,$.
\end{proof}

\pars
We leave it as an exercise to the reader to prove the following easy
lemma:
\begin{lemma}                               \label{triv}
Let $E|K$ be a finitely generated field extension and $P$ a trivial
place on $E$. Then $(\tilde{E}|K,P)$ has normal-uniformization over $E$.
If in addition $E|K$ is separable, then $(E|K,P)$ is uniformizable
and $(E\sep|K,P)$ has Galois-uniformization over $E$.
\end{lemma}

%
%
\section{Uniformizable valued field extensions}        \label{sectufe}
In this section, we will present various finitely generated valued field
extensions which are uniformizable.

%
%
\subsection{Rational function fields with Abhyankar places}
The following lemma was proved (but not explicitly stated) by Zariski in
[Z1] for subgroups of $\R$, using the algorithm of Perron. We leave it
as an easy exercise to the reader to prove the general case by induction
on the rank of the ordered abelian group. However, an instant proof of
the lemma can also be found in [EL] (Theorem 2.2).
\begin{lemma}                               \label{perron}
Let $\Gamma$ be a finitely generated ordered abelian group. Take any
non-negative elements $\alpha_1,\ldots,\alpha_\ell\in\Gamma$. Then there
exist positive elements $\gamma_1,\ldots,\gamma_{\rho}\in\Gamma$ such
that $\Gamma=\Z \gamma_1 \oplus\ldots\oplus\Z\gamma_{\rho}$ and every
$\alpha_i$ can be written as a sum $\sum_{j}^{} n_{ij}\gamma_j$ with
non-negative integers $n_{ij}\,$.
\end{lemma}

The foregoing lemma and Theorem~\ref{prelBour} are the main ingredients
in the proof of the next proposition.
We consider a function field $F|K$ and a place $P$ of $F$ such that
$v_P^{ }K$ is a convex subgroup of $v_P^{ }F$. The latter always holds
if $P$ is trivial on $K$ since then, $v_P^{ }K=\{0\}$. We take elements
$x_1,\ldots,x_{\rho}$ in $F$ such that $v_P^{ }x_1,\ldots, v_P^{ }
x_{\rho}$ form a maximal set of rationally independent elements in
$v_P^{ } F$ modulo $v_P^{ }K$. Further, we take elements $y_1,\ldots,
y_{\tau}$ in $F$ such that $y_1P,\ldots, y_{\tau} P$ form a
transcendence basis of $FP$ over $K$.
\begin{proposition}                         \label{prep}
In the situation described above, $(K(x_1,\ldots,x_{\rho},y_1,\ldots,
y_{\tau})|K,P)$ is uniformizable. More precisely, the
transcendence basis $T=\{t_1,\ldots,t_s\}$ can be chosen of the form
$\{x'_1,\ldots,x'_{\rho}, y_1,\ldots, y_{\tau}\}$, where $x'_1,\ldots,
x'_{\rho}\in {\cal O}_{K(x_1,\ldots,x_{\rho})}$ and for some $c\in
{\cal O}_K$ (with $c=1$ if $P$ is trivial on $K$), the elements
$cx'_1,\ldots,cx'_{\rho}$ generate the same multiplicative subgroup of
$K(x_1,\ldots,x_{\rho})^\times$ as $x_1, \ldots,x_{\rho}\,$. If $c\in
{\cal O}_K$ such that $vc'\geq vc$, then $c$ can be replaced by $c'$.
\end{proposition}
\begin{proof}
For the proof, we set $F=K(x_1,\ldots,x_{\rho},y_1,\ldots,y_{\tau})$.
By Theorem~\ref{prelBour}, we know that
\[v_P^{ }F\>=\>\Z v_P^{ }x_1\oplus\ldots\oplus\Z v_P^{ }x_{\rho}
\oplus v_P^{ }K\;.\]
Thus,
\[v_P^{ }F/v_P^{ }K\ni \mu_1(v_P^{ }x_1+v_P^{ }K)+\ldots+\mu_{\rho}
(v_P^{ }x_{\rho}+v_P^{ }K)\>\mapsto\> x_1^{\mu_1}\cdot\ldots\cdot
x_{\rho}^{\mu_{\rho}}\hspace{1cm}  (\mu_1,\ldots,\mu_{\rho}\in\Z)\]
is an isomorphism from $v_P^{ }F/v_P^{ }K$ onto the multiplicative
subgroup of $F^\times$ generated by $x_1,\ldots,x_{\rho}\,$. We denote
this group by ${\cal X}$.

\pars
Now let $\zeta_1,\ldots,\zeta_m\in {\cal O}_F\,$. Take $\zeta$ to
be any of these elements and write $\zeta=f/g$ with polynomials $f,g\in
K[x_1,\ldots,x_{\rho}, y_1,\ldots,y_{\tau}]$. Write
\[f\>=\>\sum_{i=1}^{d} c_i\, x_1^{\mu_{1,i}}\cdot\ldots\cdot
x_{\rho}^{\mu_{\rho,i}}\cdot y_1^{\nu_{1,i}}\cdot\ldots\cdot
y_{\tau}^{\nu_{\tau,i}} \;\;\mbox{ and }\;\;
g\>=\>\sum_{i=1}^{d'} c'_i\, x_1^{\mu'_{1,i}}\cdot\ldots\cdot
x_{\rho}^{\mu'_{\rho,i}}\cdot y_1^{\nu'_{1,i}}\cdot\ldots\cdot
y_{\tau}^{\nu'_{\tau,i}}\]
as sums of monomials (in such a way that in either polynomial, two
different monomials differ in at least one exponent). Then by
Theorem~\ref{prelBour}, the value of $f$ is equal to the least of the
values of its monomials, say, to the one of the first. Similarly for
$g$. So we can write
\begin{equation}                            \label{preper}
\zeta\>=\>
%
\frac{\sum_{i}\frac{c_i}{c'_1}\,x_1^{\mu_{1,i}-\mu'_{1,1}}\cdot\ldots
\cdot x_{\rho}^{\mu_{\rho,i}-\mu'_{\rho,1}}\cdot y_1^{\nu_{1,i}}\cdot
\ldots\cdot y_{\tau}^{\nu_{\tau,i}}}{\sum_{i}^{} \frac{c'_i}{c'_1}\,
x_1^{\mu'_{1,i}-\mu'_{1,1}}
\cdot\ldots\cdot x_{\rho}^{\mu'_{\rho,i}-\mu'_{\rho,1}}\cdot
y_1^{\nu'_{1,i}}\cdot\ldots\cdot y_{\tau}^{\nu'_{\tau,i}}}
\end{equation}
where the denominator has value 0 and the summands appearing in it
all have value $\geq 0$. Since $\zeta\in {\cal O}_F\,$, also the
numerator has value $\geq 0$, and the same must thus be true for all
its summands.
The only obstruction is that some of the $x_i$'s may appear with
negative exponents in some summands, and that, if $P$ is non-trivial on
$K$, some of the $c_i/c'_1$ or $c'_i/c'_1$ may have negative value.

\pars
We collect all summands of the form $h=cx_1^{\mu_1}\cdot\ldots\cdot
x_{\rho}^{\mu_{\rho}}\cdot y_1^{\nu_1}\cdot\ldots\cdot
y_{\tau}^{\nu_{\tau}}$ which appear in the numerator or denominator of
(\ref{preper}), and all products of the form $\xi=x_1^{\mu_1}\cdot
\ldots\cdot x_{\rho}^{\mu_{\rho}}\in {\cal X}$ which appear in these
summands. We do the same for all elements $\zeta_1, \ldots,\zeta_m\,$.
In this way, we obtain finitely many elements $h_1,\ldots,h_k$, all of
them having non-negative value, and corresponding elements
$\xi_1,\ldots,\xi_k\in {\cal X}$. We note that
$v_P^{ }\xi_j\in v_P^{ }h_j+v_P^{ }K$, and as $v_P^{ }K$ is a convex
subgroup of $v_P^{ }F$ and $v_P^{ }h_j$ is non-negative, it follows that
$v_P^{ } \xi_j+v_P^{ }K$ is non-negative in $v_P^{ }F/v_P^{ }K$, for the
induced order. After adding some suitably chosen elements of the form
$x_i$ or $x_i^{-1}$ (depending on whether $v_P^{ }x_i$ or $v_P^{ }
x_i^{-1}$ is positive) to the $\xi$'s if necessary, we obtain elements
$\xi_1,\ldots,\xi_\ell$ which also generate ${\cal X}$. At this point,
we apply Lemma~\ref{perron} to the non-negative values $v_P^{ }\xi_1+
v_P^{ } K,\ldots,v_P^{ }\xi_\ell+v_P^{ }K\in v_P^{ }F/v_P^{ }K$. Pulling
the result back through the above isomorphism, we find generators
$x''_1,\ldots,x''_{\rho}$ of ${\cal X}$ for which the values $v_P^{ }
x''_j +v_P^{ }K$ in $v_P^{ }F/v_P^{ }K$ are positive, and such that
every of the $\xi$'s can be written as a (unique) product of the $x''_j$
{\it with non-negative exponents}.

Let $\alpha\in v_P^{ }K$ be the minimum of the values of the
coefficients $c_1,\ldots,c_k$ appearing in the monomials $h_1,\ldots,
h_k\,$. We take $c\in {\cal O}_K$ such that $v_P^{ }c\geq -\min
\{0,\alpha\}$. If $\alpha\geq 0$, which in particular is the case if $P$
is trivial on $K$, then we can take $c:=1$. Now we set $x'_j:=
c^{-1}x''_j\,$. It follows that $h_1, \ldots,h_k\in {\cal O}_K[x'_1,
\ldots,x'_{\rho},y_1,\ldots,y_{\tau}]$. Since the values $v_P^{ } x''_j
+v_P^{ }K$ in $v_P^{ }F/v_P^{ }K$ were positive, all values $v_P^{ }
x'_j$ are positive. This remains true if $c$ is replaced by any $c'
\in {\cal O}_K$ such that $vc'\geq vc$.

\pars
As we can read off from (\ref{preper}), every $\zeta_j$ can be written
as $\zeta'_j/\zeta''_j\,$, where $\zeta'_j,\zeta''_j$ lie in ${\cal O}_K
[x'_1,\ldots,x'_{\rho},y_1,\ldots,y_{\tau}]\,$, with $v_P^{ }\zeta''_j
=0$. We set $s=\rho+\tau$, $t_i:=x'_i$ for $1\leq i\leq\rho$ and
$t_{\rho+i}:=y_i$ for $1\leq i\leq\tau$. Next, we set $n:=m$ and put
$\eta_j:= \zeta_j$ and $f_j(X_1,\ldots,X_n):= \zeta''_jX_j-\zeta'_j\,$,
and we are done.
\end{proof}

%
%
\subsection{Immediate simple transcendental extensions}

\begin{lemma}                               \label{ist}
Let $(K(z)|K,P)$ be an immediate transcendental extension. If
(\ref{trat}) holds, then $(K(z)|K,P)$ is uniformizable.
\end{lemma}
\begin{proof}
Let $\zeta_1,\ldots,\zeta_m\in {\cal O}_{K(z)}\,$ and write $\zeta_j=
f_j(z)/g_j(z)$ with polynomials $f_j(z)$, $g_j(z)\in K[z]$. We apply
Lemma~\ref{lvpol} to these finitely many polynomials and choose
$\tilde{z}=\frac{z-a}{b}$ according to this lemma. Then by (\ref{vpol}),
for every $j$ we can find $i_j,k_j$ such that $v_P^{ }f_j(z)= v_P^{ }
f_j^{[i_j]} (a)\,b^{i_j} =\min_i v_P^{ } f_j\iT (a)\,b^i$
and $v_P^{ }g_j(z)= v_P^{ }g_j^{[k_j]} (a)\,b^{k_j}=
\min_i v_P^{ }g_j\iT (a)\,b^i$. Thus, we can write
\[\zeta_j\>=\>\frac{f_j^{[i_j]}(a)\,b^{i_j}} {g_j^{[k_j]}(a)\,b^{k_j}}
\cdot \frac{\tilde{f}_j(\tilde{z})}{\tilde{g}_j(\tilde{z})}\]
where $\tilde{f}_j,\tilde{g}_j$ are polynomials with coefficients in
${\cal O}_K$ and $v_P^{ }\tilde{f}_j(\tilde{z})=0=v_P^{ } \tilde{g}_j
(\tilde{z})$. Note that also the first fraction is an element of
${\cal O}_K$ since its value is equal to $v_P^{ }\zeta_j$ and $\zeta_j
\in {\cal O}_{K(z)}$ by assumption. Now we set $t_1:=\tilde{z}$,
$n:=m$, $\eta_j:=\zeta_j$ and
\[f_j(X_1,\ldots,X_n)\>:=\> \tilde{g}_j(t_1)\, X_j \,-\,
\frac{f_j^{[i_j]}(a)\,b^{i_j}}{g_j^{[k_j]}
(a)\,b^{k_j}}\cdot\tilde{f}_j(t_1)\]
for $1\leq j\leq m$, and we are done.
\end{proof}

%
%
\subsection{Extensions within the completion}      \label{sectec}
\begin{lemma}                               \label{comp}
Every finite separable-algebraic extension of a valued field within its
completion is uniformizable.
\end{lemma}
\begin{proof}
Take any separable-algebraic extension $(L|K,P)$ such that $(L,P)$ lies
in the completion of $(K,P)$. Further, take $\zeta_1,\ldots,\zeta_m\in
{\cal O}_{L}\,$. Let $\zeta$ be any of these elements.
We extend $v_P^{ }$ to the algebraic closure of the completion.
Since $L|K$ is separable, we can write the minimal
polynomial of $\zeta$ in the form $(X-\zeta)(X-\sigma_1\zeta)\ldots
(X-\sigma_k\zeta)$ where $\sigma_j$ are automorphisms in
Aut$(K\sep|K)$ and all $\sigma_j\zeta$ are distinct from $\zeta$.
Since $\zeta$ lies in the completion of $(K,P)$, it follows that there
is some $a\in K$ such that $v_P^{ }(\zeta -a)>v_P^{ }(\sigma_j\zeta-a)$
for all $j$, and that $v_P^{ }(\zeta -a)\geq 0$. Since $\zeta\in
{\cal O}_L\,$, the latter implies that $a\in {\cal O}_K\,$. On
the other hand, $v_P^{ }(\zeta -a)\in v_P^{ }K(\zeta)=v_P^{ }K$ and
$K(\zeta)P=KP$, so we can choose $b\in {\cal O}_K$ such that
$v_P^{ } (\frac{\zeta-a}{b}-1) >0$. Thus, $v_P^{ }\frac{b}{\zeta-a}
=0$ and $v_P^{ }\sigma_j\frac{b} {\zeta-a} = v_P^{ }
\frac{b}{\sigma_j\zeta-a}>0$. Therefore, the reduction $hP$ of the
minimal polynomial $h$ of $\frac{b}{\zeta-a}$ over
$K$ is the polynomial $X^{k+1}-X^k$ which has $1=\frac{b}{\zeta-a}P$ as
a simple root. We set $n:=3m$. If $\zeta$ was $\zeta_j$ then we set
$\eta_{3j-2}:= \frac{b}{\zeta-a}$, $f_{3j-2}\>:=\>h(X_{3j-2})$,
$\eta_{3j-1}:= \frac{\zeta-a}{b}$, $f_{3j-1}\>:=\>X_{3j-2}X_{3j-1}-1$,
$\eta_{3j}:= \zeta=b\frac{\zeta-a}{b}+a\in {\cal O}_K[\frac{\zeta-a}
{b}]$, $f_{3j}\>:=\>X_{3j}-bX_{3j-1}-a$, and we are done.
\end{proof}

We can drop the condition that the extension be algebraic:
\begin{proposition}                               \label{comp1}
Every finitely generated separable extension of a valued field within
its completion is uniformizable.
\end{proposition}
\begin{proof}
It suffices to prove the assertion for every finitely generated
separable extension $(L|K,P)$ within the completion of $(K,P)$. As $L|K$
is finitely generated and separable, we can choose a transcendence basis
$z_1,\ldots,z_n$ such that $L|K(z_1,\ldots,z_n)$ is separable-algebraic.
By induction on the transcendence degree, using Lemma~\ref{close},
Lemma~\ref{ist} and transitivity (Lemma~\ref{basictrans}), we find that
$(K(z_1,\ldots,z_n)|K,P)$ is uniformizable. By
the foregoing lemma, the same holds for $(L|K(z_1,\ldots,z_n),P)$. Now
our assertion follows by transitivity.
\end{proof}

%
%
\subsection{Extensions within the henselization}      \label{secteh}
The henselization of a valued field $(K,P)$ is always a
separable-algebraic extension. If $(K,P)$ has rank 1, then moreover, the
henselization lies in the completion of $(K,P)$ (since in this case the
completion is henselian, cf.\ [R], [K2]). Therefore, Lemma~\ref{comp}
yields:
\begin{corollary}                               \label{hens1}
Assume that $(K,P)$ has rank 1. Then every finite extension of $(K,v)$
within its henselization is uniformizable.
\end{corollary}

We give a typical application:
\begin{corollary}                           \label{Abhy1}
If $P$ is a rational Abhyankar place of rank 1 of a function field
$F|K$, then $(F|K,P)$ is uniformizable.
\end{corollary}
\begin{proof}
By Theorem~\ref{hrwtd}, $(F|K,P)$ is henselian generated and there are
$x_1,\ldots,x_{\rho}\in F$ as in the assertion of that theorem such that
$(F,P) \subset (K(x_1,\ldots,x_{\rho}),P)^h$. By the foregoing corollary
it follows that $(F|K(x_1,\ldots,x_{\rho}),P)$ is uniformizable. By
Proposition~\ref{prep}, $(K(x_1,\ldots,x_{\rho})|K,P)$ is uniformizable.
Now our assertion follows by transitivity.
\end{proof}

\parm
To treat the case of a rank higher than 1, we use a
well known lemma about composite places (cf.\ [R] or [K2]).
\begin{lemma}                              \label{hensPQ}
Suppose that the place $P$ of $K$ is composite: $P=Q\ovl{Q}\,$. Then
$(K,P)$ is henselian if and only if $(K,Q)$ and $(KQ,\ovl{Q})$ are.
If $(KQ,\ovl{Q})$ is henselian, then the henselization of $K$ with
respect to $P$ is equal to the henselization of $K$ with respect to
$Q\,$ (as fields).
\end{lemma}

If in this situation, $Q$ has rank 1, the henselization of $K$ with
respect to $Q$ lies in its completion with respect to $Q$. Since $P=Q
\ovl{Q}$, it follows from general valuation theory that the completion
of $K$ with respect to $Q$ is equal to the completion of $K$ with
respect to $P$. So the henselization of $K$ with respect to $P$ lies in
the completion of $K$ with respect to $P$. Hence, we obtain the
following corollary from Lemma~\ref{comp}:
\begin{corollary}                           \label{corPQ}
Assume that $P$ is a place of $K$ and $P=Q\ovl{Q}$ such that $(K,Q)$ has
rank 1 and $(KQ,\ovl{Q})$ is henselian. Then every finite extension of
$(K,P)$ within its henselization is uniformizable.
%
\end{corollary}

%
%
\subsection{Immediate extensions}
\begin{proposition}                         \label{isofstF}
Take a separable-algebraically closed field $K$, a separable function
field $F|K$ of transcendence degree 1, and a place $P$ on $F$ of rank 1
such that $(F|K,P)$ is an immediate extension. Then $(F|K,P)$
is uniformizable.

The assertion also holds if $P=Q\ovl{Q}$ such that $(F,Q)$ has rank 1
and $FQ=KQ$.
\end{proposition}
\begin{proof}
By Theorem~\ref{stt3}, $(F|K,P)$ is henselian generated. That is,
there is some $z\in F$ such that $F\subset K(z)^h$. Since $K$ is
separable-algebraically closed, Lemma~\ref{close} shows that condition
(\ref{trat}) holds. Therefore, Lemma~\ref{ist} shows that $(K(z)|K,P)$
is uniformizable. By Corollary~\ref{hens1}, the same holds for $(F|K(z),
P)$. Hence by transitivity, $(F|K,P)$ is uniformizable.

If $P=Q\ovl{Q}$ such that $(F,Q)$ has rank 1 and $FQ=KQ$, then we employ
Corollary~\ref{corPQ} in the place of Corollary~\ref{hens1}. This is
possible since $FQ=KQ$ implies that $K(z)Q=KQ$ and thus, being equal
to the residue field of a separable-algebraically closed field, $K(z)Q$
is itself separable-algebraically closed and hence henselian under every
valuation.
\end{proof}

%
%
\section{Galois- and normal-uniformization}      \label{sectGnu}
In this section, we will present valued field extensions which admit
Galois-uniformization or normal-uniformization.

%
%
\subsection{Abhyankar places of rank 1}
\begin{proposition}                         \label{r1zdAbh}
Take a function field $E|K$ and a zero-dimensional Abhyankar place $P$
of $E|K$ of rank 1. Then $(\tilde{E}|K,P)$ has normal-uniformization
over $E$. If $K$ is perfect, then $(E\sep|K,P)$ has
Galois-uniformization over $E$.
\end{proposition}
\begin{proof}
For given $\zeta_1,\ldots,\zeta_m\in {\cal O}_{E\sep}$, we take $F$ to
be the normal hull of $E(\zeta_1,\ldots,\zeta_m)$ over $E$; then $F|E$
is a finite Galois extension. Since $K$ is assumed to be perfect and
$EP|K$ to be algebraic, $FP|K$ is a separable-algebraic extension. By
Lemma~\ref{ifcont}, for $F':=F.K\sep$ we have that $F'P= (FP)\sep=
K\sep$. Thus, Corollary~\ref{Abhy1} shows that $(F'|K\sep,P)$ is
uniformizable. This proves that $(E\sep|K\sep,P)$
has Galois-uniformization over~$E$.

On the other hand, $P$ is trivial on $K\sep$. Hence, Lemma~\ref{triv}
tells us that $(K\sep|K,P)$ has Galois-uniformization over $K$.
Now our assertion follows by transitivity.

The proof for normal-uniformization is similar.
\end{proof}

\begin{proposition}                         \label{AbhPQ}
Take a function field $E|K$ and a place $P$ of $E$. Assume that
$P=Q\ovl{Q}$ such that $Q$ is a zero-dimensional Abhyankar place
of $E|K$ of rank 1. Then $(\tilde{E}|\tilde{K},P)$ has
normal-uniformization over $E$.
\end{proposition}
\begin{proof}
For given $\zeta_1,\ldots,\zeta_m\in {\cal O}_{\tilde{E}}$, we take $F$
to be the normal hull of $E(\zeta_1,\ldots,\zeta_m)$ over $E$. We set
$F':=F.\tilde{K}$. As $Q$ is zero-dimensional, we obtain that $F'Q=
\tilde{K}$. Hence, $Q$ is a rational Abhyankar place of $F'|\tilde{K}$.
By Theorem~\ref{hrwtd}, $(F'|\tilde{K},Q)$ is henselian generated: there
are $x_1^Q,\ldots,x_k^Q\in E$ such that $v_Q^{ }x_1^Q,\ldots, v_Q^{ }
x_k^Q$ is a maximal set of rationally independent values in $v_Q^{ }F'$,
and $(F',Q)$ is contained in the henselization of $(F_0,Q)$, where
$F_0:=\tilde{K}(x_1^Q,\ldots,x_k^Q)$. But $F_0Q=\tilde{K}$ by
Theorem~\ref{prelBour}, and $(\tilde{K},\ovl{Q})$ is henselian. Hence by
Corollary~\ref{corPQ}, $(F'|F_0,P)$ is uniformizable.

\pars
Since $Q$ is trivial on $\tilde{K}$ and since $F'Q=\tilde{K}$, we know
that $v_P^{ }\tilde{K} = v_{\ovl{Q}}^{ }\tilde{K} = v_{\ovl{Q}}^{ }
(F'Q)$ is a convex subgroup of $v_P^{ }F'$, and that $v_Q^{ }F'=v_P^{ }
F'/v_P^{ }\tilde{K}$. Consequently, our choice of the $x_i^Q$'s yields
that $v_P^{ }x_1^Q,\ldots, v_P^{ } x_k^Q$ form a maximal set of
rationally independent elements in $v_P^{ } F'$ modulo $v_P^{ }
\tilde{K}$. Hence by Proposition~\ref{prep}, $(F_0|\tilde{K},P)$
is uniformizable. By transitivity, the same holds for
$(F'|\tilde{K},P)$. This proves our assertion.
\end{proof}

\begin{proposition}                           \label{corr1nuAbh}
Take a function field $E|K$ and an Abhyankar place $P$ of $E|K$ of
rank~1. Then $(\tilde{E}|K,P)$ has normal-uniformization over $E$.
\end{proposition}
\begin{proof}
We choose $y_1,\ldots,y_{\tau}\in E$ such that $y_1P,\ldots,y_{\tau}P$
is a transcendence basis of $EP|K$. We take $K'$ to be the algebraic
closure of $K(y_1,\ldots,y_{\tau})$. We extend $P$ to $\tilde{E}$. Then
$P$ induces an isomorphism on $K'$. Passing to an equivalent place if
necessary, we can assume that $\tilde{E}P=K'$. Hence by
Proposition~\ref{r1zdAbh}, $(\tilde{E}|K',P)$ has normal-uniformization
over $E$. By Lemma~\ref{triv}, $(K'|K,P)$ has normal-uniformization over
$K(y_1,\ldots,y_{\tau})$. Now our assertion follows by transitivity.
\end{proof}

%
%
\subsection{Immediate extensions}
\begin{proposition}                         \label{isofst}
Take a function field $E|K$ and a place $P$ on $\tilde{E}$ of rank 1
such that $v_P^{ }E/v_P^{ }K$ is a torsion group and $EP|KP$ is
algebraic. Then the immediate extension $(\tilde{E}|\tilde{K},P)$ has
normal-uniformization over $E$, and if $E|K$ is separable, then the
immediate extension $(E\sep|K\sep,P)$ has Galois-uniformization over
$E$.

These assertions also hold if $P=Q\ovl{Q}$ such that $(E,Q)$ has rank 1,
$v_Q^{ }E/v_Q^{ }K$ is a torsion group and $EQ|KQ$ is algebraic.
\end{proposition}
\begin{proof}
We give the proof for Galois-uniformization. We proceed by induction on
the transcendence degree. The case of transcendence degree 1 is covered
by Proposition~\ref{isofstF}: For given $\zeta_1,\ldots,\zeta_m\in
{\cal O}_{E\sep}\,$, we take $F$ to be the normal hull of $E(\zeta_1,
\ldots, \zeta_m)$ over $E$, which is a finite Galois extension of $E$.
Then we apply Proposition~\ref{isofstF} to $(F.K\sep|K\sep,P)$. We
observe that since $v_P^{ }E/v_P^{ }K$ is a torsion group and $EP|KP$ is
algebraic by hypothesis, Lemma~\ref{KsacP} implies that the extension
$(E\sep|K\sep,P)$ and hence also its subextension $(F.K\sep|K\sep,P)$
are immediate. For the case of $P=Q\ovl{Q}$ with $Q$ non-trivial, it
also implies that $K\sep Q=\widetilde{KQ}$ and $E\sep Q=\widetilde{EQ}$.
Hence, our assumption that $EQ|KQ$ is algebraic implies that $E\sep
Q=K\sep Q$.

So let us now assume that $\trdeg E|K=n>1$ and that our assertion is
true for transcendence degree $<n$. We take a separating transcendence
basis $T$ of $E|K$. Then we pick a subset $T_0\subset T$ such that
$\trdeg E|E_0=1$ for $E_0:=K(T_0)\subset E$. It follows that
$E.E_0\sep|E_0\sep$ is a separable function field of transcendence
degree $1$ and that $E_0|K$ is a separable function field of
transcendence degree $n-1$. As $(E_0|K,P)$ is a subextension of
$(E|K,P)$, $v_P^{ }E_0/v_P^{ }K$ is a torsion group, $E_0P|KP$ is
algebraic, and $(E_0,P)$ will have rank at most 1 if $(E,P)$ has rank 1.
But the fact that $v_P^{ }E/v_P^{ }K$ is a torsion group implies that
$(K,P)$ has the same rank as $(E,P)$. This shows that $(E_0,P)$ has rank
1 if $(E,P)$ has rank 1. Similarly, if $P=Q\ovl{Q}$ with $EQ=KQ$ and
$(E,Q)$ has rank 1, then $E_0Q=KQ$ and $(E_0,Q)$ will have rank at most
1. But the fact that $v_P^{ }E/v_P^{ }K$ is a torsion group also implies
that $(K,Q)$ has the same rank as $(E,Q)$. Hence in this case, $(E_0,Q)$
has rank 1 if $(E,Q)$ has rank 1. We have shown that also $(E_0|K,P)$
satisfies the assumptions of our proposition. As $(E.E_0\sep)\sep
=E\sep$, our induction hypothesis yields that $(E\sep|E_0\sep,P)$ has
Galois-uniformization over $E$ and that $(E_0\sep|K\sep,P)$ has
Galois-uniformization over $E_0\,$. Hence by transitivity, $(E\sep
|K\sep,P)$ has Galois-uniformization over $E$.

The proof for normal-uniformization is similar: instead of
separable-algebraic closures we use algebraic closures.
\end{proof}

%
%
\subsection{Places of rank 1}
\begin{proposition}                         \label{r1zd}
Take a function field $E|K$ and a zero-dimensional place $P$ of $E|K$ of
rank 1. If $K$ is perfect, then $(E\sep|K,P)$ has Galois-uniformization
over $E$.
\end{proposition}
\begin{proof}
By Lemma~\ref{tb}, we can choose a separating transcendence basis of
$E|K$ which contains elements $x_1,\ldots,x_{\rho}\in E$ such that
$v_P^{ } x_1,\ldots,v_P^{ } x_{\rho}$ is a maximal set of rationally
independent elements in $v_P^{ }E\,$. We set $E_0:=K(x_1,\ldots,
x_{\rho})$. Then $E|E_0$ is separable. Further, $v_P^{ }E/v_P^{ }E_0$ is
a torsion group by Theorem~\ref{prelBour}. By the same theorem,
$E_0P=K$. By assumption, $EP$ is algebraic over $K=E_0P$. Hence by
Proposition~\ref{isofst}, the extension $(E\sep |E_0\sep,P)$ has
Galois-uniformization over $E$. By Proposition~\ref{r1zdAbh},
$(E_0\sep|K,P)$ has Galois-uniformization over $E_0\,$. Now our
assertion follows by transitivity.
\end{proof}


%
%
\subsection{Places of arbitrary rank}
\begin{proposition}                         \label{rr1}
Take a function field $E|K$ and a place $P$ on $E$. Assume that
$P=Q\ovl{Q}$ such that $Q$ is a place of $E|K$ of rank 1. Then
$(\tilde{E}|\tilde{K},P)$ has normal-uniformization over~$E$.
\end{proposition}
\begin{proof}
We choose $x_1^Q,\ldots,x_k^Q\in E$ such that $v_Q^{ }x_1^Q,\ldots,
v_Q^{ }x_k^Q$ is a maximal set of rationally independent values in
$v_Q^{ }E$. We set $L:= K(x_1^Q,\ldots,x_k^Q)\subseteq E$. Then by
Theorem~\ref{prelBour}, $LQ=KQ=K$ and $v_Q^{ }E/v_Q^{ }L$ is a torsion
group. Hence by Proposition~\ref{isofst}, $(\tilde{E}|\tilde{L},P)$ has
normal-uniformization over $E$.

By construction, $(L|K,P)$ and $Q$ satisfy the assumptions of
Proposition~\ref{AbhPQ}. This yields that $(\tilde{L}|\tilde{K},P)$
has normal-uniformization over $L$. Now our assertion follows by
transitivity.
\end{proof}

\begin{proposition}                         \label{nur>}
Take a function field $E|K$ and a place $P$ of $E|K$. Then
$(\tilde{E}|K,P)$ has normal-uniformization over $E$.
\end{proposition}
\begin{proof}
Since the rank of $(E,P)$ is finite (cf.\ Section~\ref{sectrank}), we
can proceed by induction on this rank. Assume that $P=Q\ovl{Q}$ such
that $Q$ is a place of $E|K$ of rank 1, with $\ovl{Q}$ possibly trivial.
We take $y_1^Q,\ldots,y_\ell^Q \in E$ such that $y_1^Q Q,\ldots,y_\ell^Q
Q$ is a transcendence basis of $EQ|K$. Then we set $E_1:=K(y_1^Q,
\ldots,y_\ell^Q)\subseteq E$ and $E':=E.\tilde{E}_1\subset\tilde{E}$.
Since $E'|E$ is algebraic, so is $E'Q|EQ$. On the other hand, $EQ$ is
algebraic over $E_1Q$ by construction, and therefore, $\tilde{E}_1 Q=
\widetilde{E_1 Q}$ is equal to the algebraic closure of $EQ$. As
$\tilde{E}_1 Q\subseteq E'Q$, this shows that $E'Q=\tilde{E}_1 Q$. Since
$Q$ is the identity on $K$ and $y_1^Q Q,\ldots,y_\ell^Q Q$ are
algebraically independent over $K$, it induces an isomorphism on $E_1$
and hence also on $\tilde{E}_1\,$. Passing to an equivalent place if
necessary, we can assume that $Q$ is a place of $E'|\tilde{E}_1\,$.

Since $\widetilde{E'}=\tilde{E}$, Proposition~\ref{rr1} now shows that
$(\tilde{E}|\tilde{E}_1,P)$ has normal-uniformization over $E'$, and
hence also over $E$. As the rank of $(E_1,P)$ is equal to the rank of
$(EQ,\ovl{Q})$ and thus smaller than the rank of $(E,P)$, our induction
hypothesis (or Lemma~\ref{triv}, if $\ovl{Q}$ is trivial) yields that
$(\tilde{E}_1|K,P)$ has normal-uniformization over $E_1\,$. Now our
assertion follows by transitivity.
\end{proof}

%
%
\section{Proof of the main theorems for rank 1}     \label{sectprMT1}
$\bullet$ \ {\bf Proof of Theorems~\ref{MT1}, \ref{MT1norm} for rank 1,
and of Theorem~\ref{MT4}}
\sn
Theorems~\ref{MT1} and \ref{MT1norm} can be proved by a direct
application of Proposition~\ref{nur>}. But we will use a different
approach which at the same time proves Theorem~\ref{MT4}.

Let $F|K$ be a function field and $E|K$ a subextension of the same
transcendence degree; consequently, $F|E$ is finite. Further, take a
rank 1 place $P$ of $F|K$ and any elements $\zeta_1,\ldots,\zeta_m\in
{\cal O}_F\,$. After extending this list if necessary, we can assume
that it includes generators of $F|E$. Finally, take $E_0$ to be any
Abhyankar field of $(E|K,P)$.

By Proposition~\ref{isofst}, $(\tilde{E}|\tilde{E}_0,P)$ has
normal-uniformization over $E$. Therefore, there is a finite normal
extension ${\cal F}'$ of $E$ such that ${\cal F}'.\tilde{E}_0$
contains the $\zeta$'s and $({\cal F}'.\tilde{E}_0|\tilde{E}_0,P)$ is
uniformizable with respect to the $\zeta$'s. By part b) of
Lemma~\ref{going}, there is a finite extension ${\cal E}_0|E_0$ such
that $({\cal F}'.{\cal E}_0|{\cal E}_0,P)$ is uniformizable with
respect to the $\zeta$'s.

We choose $\zeta'_1,\ldots,\zeta'_{m'}\in {\cal O}_{{\cal E}_0}$ to
consist of generators of ${\cal E}_0$ over $E_0$ and of the
uniformization coefficients of ${\cal F}'.{\cal E}_0,\zeta_1,\ldots,
\zeta_m$ in ${\cal E}_0\,$. By Proposition~\ref{corr1nuAbh},
$(\tilde{E}_0|K,P)$ has normal-uniformization over $E_0\,$. Hence there
is a finite normal extension ${\cal F}_0$ of $E_0$ and a purely
inseparable subextension ${\cal K}|K$ of ${\cal F}_0|K$ such that
${\cal F}_0$ contains $\zeta'_1,\ldots,\zeta'_{m'}$ (and hence also
${\cal E}_0$) and $({\cal F}_0|{\cal K},P)$ is uniformizable with
respect to $\zeta'_1,\ldots,\zeta'_{m'}\,$, the uniformizing
transcendence basis being a transcendence basis of ${\cal F}_0|K$. We
have that $({\cal F}'.{\cal E}_0).{\cal F}_0={\cal F}'.{\cal F}_0\,$.

Now by Lemma~\ref{basictrans}, $({\cal F}'.{\cal F}_0|{\cal K},P)$ is
uniformizable with respect to $\zeta_1,\ldots,\zeta_m\,$, with
uniformizing transcendence basis containing a transcendence basis of
${\cal F}_0|K$. With $\rho$ and $\tau$ as defined preceding to
Theorem~\ref{MT4}, the latter transcendence basis has $\rho+\tau$
many elements, and they are algebraic over $E_0\,$. We set ${\cal
F}:={\cal F}'.{\cal F}_0\,$. As ${\cal F}'|E$ is finite and normal and
${\cal F}_0|E_0$ is finite and normal, ${\cal F} |E$ is also finite and
normal. By our additional assumption on the $\zeta$'s, $F\subset
{\cal F}$. Hence also ${\cal F}|F$ is finite and normal. This proves
Theorems~\ref{MT1} and~\ref{MT1norm} in the rank 1 case, and
Theorem~\ref{MT4}.

\bn
$\bullet$ \ {\bf Proof of Theorem~\ref{MTr1E}}
\sn
If $K$ is not perfect, we will have to pass to its perfect hull to apply
Propositions~\ref{r1zd}. This perfect hull will be denoted by
$K^{1/p^\infty}$.

For the proof of Theorem~\ref{MTr1E}, we assume in addition to the
above assumptions that $P$ is zero-dimensional and that $F|E$ is
separable-algebraic, i.e., $F\subset E\sep$. We set $E_1:=
E.K^{1/p^\infty}$; then $E_1\sep=E\sep .K^{1/p^\infty}$. As $P$
is still a zero-dimensional place of $E_1\sep|K^{1/p^\infty}$,
Proposition~\ref{r1zd} shows that $(E_1\sep|K^{1/p^\infty}, P)$ has
Galois-uniformization over $E_1\,$. This gives us a finite Galois
extension ${\cal F}''$ of $E_1$ such that $({\cal F}''|
K^{1/p^\infty},P)$ is uniformizable with respect to the $\zeta$'s. Now
we take ${\cal F}'$ to be the maximal separable subextension of
${\cal F}''|E$. Then ${\cal F}'|E$ is a Galois extension and ${\cal F}''
={\cal F}'.K^{1/p^\infty}$. Thus, $({\cal F}'. K^{1/p^\infty}|
K^{1/p^\infty},P)$ is uniformizable with respect to the $\zeta$'s, and
by part b) of Lemma~\ref{going}, there is a finite subextension
${\cal K} |K$ of $K^{1/p^\infty}|K$ such that $({\cal F}'.{\cal K}|
{\cal K},P)$ is uniformizable with respect to the $\zeta$'s. We set
${\cal F}:={\cal F}'.{\cal K}$. As ${\cal F}'|E$ is a Galois extension,
so is ${\cal F}|E.{\cal K}\,$. By our additional assumption on the
$\zeta$'s, $F\subset {\cal F}$. Hence also ${\cal F}|F.{\cal K}$ is a
Galois extension. This proves Theorem~\ref{MTr1E}.         \QED

\bn
$\bullet$ \ {\bf Proof of Theorem~\ref{MTr12}}
\sn
Let $F$ be as in the hypothesis of Theorem~\ref{MT1}, with $P$ of rank
1. We set $F_1:=F.K^{1/p^\infty}$. By Lemma~\ref{tb}, we can choose a
separating transcendence basis of $F_1|K^{1/p^\infty}$ which contains
elements $x_1,\ldots,x_{\rho}$ such that $v_P^{ }x_1, \ldots,v_P^{ }
x_{\rho}$ is a maximal set of rationally independent elements in
$v_P^{ }F_1\,$, and $y_1,\ldots,y_{\tau}$ such that $y_1P,\ldots,
y_{\tau}P$ is a transcendence basis of $F_1P|K^{1/p^\infty}$. We set
$F_0:=K^{1/p^\infty}(x_1,\ldots,x_{\rho},y_1,\ldots,y_{\tau})\subseteq
F_1\,$. Then $F_1|F_0$ is separable, $v_P^{ }F_1/v_P^{ }F_0$ is a
torsion group and $F_1P|F_0P$ is algebraic.
%
%
Hence by Proposition~\ref{isofst}, the extension $(F_1\sep |F_0\sep,P)$
has Galois-uniformization over $F_1$. That is, there is a finite Galois
extension $F_2|F_1$ such that $(F_2.F_0\sep |F_0\sep,P)$ is
uniformizable with respect to the $\zeta$'s. By part c) of
Lemma~\ref{going} there is a finite Galois extension $F'_0$ of $F_0$
such that $(F_2.F'_0|F'_0,P)$ is uniformizable with respect to
$\zeta_1,\ldots,\zeta_m\,$. Then $F':=F_2.F'_0$ is a finite Galois
extension of $F_1\,$. By construction, $P$ is an Abhyankar place on
$F_0$ and $F'_0\,$.

Now we choose $y'_1,\ldots,y'_{\tau}\in F'$ such that $y'_1P,\ldots,
y'_{\tau}P$ is a separating transcendence basis of $F'_0P|
K^{1/p^\infty}$. We take $K'$ to be the separable-algebraic closure of
$K^{1/p^\infty}(y'_1,\ldots,y'_{\tau})$. Then by Lemma~\ref{ifcont}, the
place $P$ of $F'_0.K'| K'$ is rational. Hence by Corollary~\ref{Abhy1},
$(F'_0.K'|K',P)$ is uniformizable. It follows by part a)
of Lemma~\ref{going} and transitivity that $(F'.K'|K',P)$ is
uniformizable with respect to the $\zeta$'s. Hence by part c) of
Lemma~\ref{going} and transitivity, there is a finite Galois extension
$K''$ of $K^{1/p^\infty} (y'_1,\ldots,y'_{\tau})$ such that $(F'.K''|
K'',P)$ is uniformizable with respect to the $\zeta$'s. Further, $P$ is
trivial on $K''$. Hence by Lemma~\ref{triv}, $(K''|K^{1/p^\infty},P)$ is
uniformizable. By transitivity, $(F'.K''|K^{1/p^\infty} ,P)$ is
uniformizable with respect to the $\zeta$'s. Observe that $F'.K''|K'$ is
a finite Galois extension.

We denote by ${\cal F}'$ the maximal separable subextension of $F$ in
$F'$, and by ${\cal F}''$ the maximal separable subextension of
${\cal F}'$ in $F'.K''$. Then ${\cal F}'|F$ and ${\cal F}''|{\cal F}'$
are Galois extensions, and we have that $F'.K''={\cal F}''.
K^{1/p^\infty}$. By part b) of Lemma~\ref{going} there is a finite
purely inseparable extension ${\cal K}$ of $K$ such that $({\cal F}''.
{\cal K}|{\cal K},P)$ is uniformizable with respect to
$\zeta_1,\ldots,\zeta_m\,$. As ${\cal F}:={\cal F}''.{\cal K}$ is a
Galois extension of the Galois extension ${\cal F}'.{\cal K}$ of
$F.{\cal K}$, this proves our theorem.                     \QED

\n
{\bf Remark:} If we could choose $y'_1,\ldots,y'_{\tau}$ in $F_1$, then
$F'.K''|F$ would be a Galois extension. But the algebraic extension
$F'_0P|F_0P$ may well be inseparable.

\bn
$\bullet$ \ {\bf Proof of Theorem~\ref{MTr1A}}
\sn
Assume that $P$ is an Abhyankar place of $F|K$ and that $(F,P)$ has
rank 1. Let us first note that part c) of Theorem~\ref{MTr1A} follows
directly from Corollary~\ref{Abhy1}. For the remaining cases, we proceed
as follows.

We set $K_1:=K$ if $FP|K$ is separable, and $K_1:=K^{1/p^\infty}$
otherwise. Then we set $F_1:=F.K_1\,$. It follows that $F_1P|K_1$ is
separable. We choose $y_1,\ldots,y_{\tau}\in F_1$ such that $y_1P,
\ldots, y_{\tau}P$ is a separating transcendence basis of $F_1P|K_1$.
Then we take $K'$ to be the separable-algebraic closure of $K_1(y_1,
\ldots,y_{\tau})$ and set $F':=F.K'=F_1.K'$. Since $P$ is trivial on
$K$, it is trivial on $K_1$ and Theorem~\ref{prelBour} yields that it is
also trivial on $K_1(y_1,\ldots,y_{\tau})$ and thus on $K'$. By
Lemma~\ref{ifcont}, $F'P=K'$. Therefore, Corollary~\ref{Abhy1} shows
that $(F'|K',P)$ is uniformizable. Hence by part c)
of Lemma~\ref{going}, for given $\zeta_1,\ldots,\zeta_m\in
{\cal O}_F$ there is a finite Galois extension ${\cal F}'$ of
$K_1(y_1,\ldots,y_{\tau})$ such that $(F. {\cal F}'|{\cal F}',P)$
is uniformizable with respect to the $\zeta$'s.

Since $P$ is trivial on ${\cal F}'$ and ${\cal F}'|K_1$ is separable,
Lemma~\ref{triv} shows that $({\cal F}'|K_1,P)$ is uniformizable.
By transitivity, $(F.{\cal F}'|K_1,P)$ is uniformizable with
respect to the $\zeta$'s.

If $FP|K$ is separable and hence $K_1=K$ by definition, then
$K_1(y_1,\ldots,y_{\tau})\subset F$ and thus, ${\cal F}:=F.{\cal F}'$
is a finite Galois extension of $F$. This proves part a) of our theorem.

For the remaining cases, we proceed as follows. Since $F.K_1$ is a
purely inseparable extension of $F$, there is some $\nu\in\N$
such that $F_0:=K(y_1^{p^\nu},\ldots,y_{\tau}^{p^\nu})\subseteq F$.
The extension ${\cal F}'|F_0$ is algebraic. Its maximal separable
subextension ${\cal F}''|F_0$ is a finite Galois extension, and
$F.{\cal F}'=F.{\cal F}''.K_1\,$. Hence, $(F.{\cal F}''.K_1|K_1,P)$ is
uniformizable with respect to the $\zeta$'s. By part b) of
Lemma~\ref{going}, there is a finite subextension ${\cal K}|K$ of
$K_1|K$ such that $(F.{\cal F}''.{\cal K}|{\cal K},P)$ is uniformizable
with respect to the $\zeta$'s. Further, ${\cal F}:=F.{\cal F}''.
{\cal K}$ is a finite Galois extension of $F.{\cal K}$. This proves
the general assertion of the theorem. If in addition $P$ is
zero-dimensional, then there are no $y$'s and we can take $F_0=K$. In
this case, ${\cal F}''$ is a finite Galois extension of $K$, which
proves part b) of the theorem.
\QED

\bn
$\bullet$ \ {\bf Proof of Theorem~\ref{MTdisc}}
\sn
Assume that $P$ is a rational discrete place of $F|K$. We write $v_P^{ }
F=\Z$. We choose a set of generators of $F$ over $K$ such that each of
these generators has $v_P^{ }$-value $1$. By Lemma~\ref{rfsepffsep},
$F|K$ is separable and thus, we can choose a separating
transcendence basis $t_1,\ldots,t_s$ for $F|K$ from these generators
(cf.\ [L], Ch.\ X, \S6, Prop.~5). Then $(F|K(t_1),P)$ is immediate.
Since $v_P^{ }F=\Z$, this implies that $(F,P)$ lies in the completion
of $(K(t_1),P)$. Hence by Proposition~\ref{comp1}, the separable
extension $(F|K(t_1),P)$ is uniformizable. By Proposition~\ref{prep},
also $(K(t_1)|K,P)$ is uniformizable. Hence by transitivity, $(F|K,P)$
is uniformizable. This proves Theorem~\ref{MTdisc}.             \QED

%
%
\section{Proof of the main theorems for arbitrary
rank}                                          \label{sectprMT}
$\bullet$ \ {\bf Proof of Theorems~\ref{MT1} and~\ref{MT1norm} for
arbitrary rank}
\sn
The proof is a direct application of Proposition~\ref{nur>}.

\bn
$\bullet$ \ {\bf Proof of Theorem~\ref{MTr>AG}}
\sn
Assume that $P$ is an Abhyankar place of rank $r>1$ of the function
field $F|K$. Since the rank of $(F|K,P)$ is finite (cf.\
Section~\ref{sectrank}), we can proceed by induction on the rank. We
take a maximal proper convex subgroup $H$ of $v_P^{ }K$. Then $v_P^{ }
K/H$ is archimedean ordered. We write $P=Q\ovl{Q}$, where $Q$ is a place
of $F|K$ of rank 1 with value group $v_P^{ }K/H$, and $\ovl{Q}$ is a
place on $FQ$ with value group $H$. By Lemma~\ref{PQAbh}, $Q$ and
$\ovl{Q}$ are Abhyankar places. Hence by Corollary~\ref{fingentb},
$v_P^{ }K$ and $H$ are finitely generated.
Now we choose $x_1,\ldots,x_k$ such that
$v_P^{ } K/H =\Z (v_P^{ }x_1 +H)\oplus \ldots\oplus\Z (v_P^{ }x_k+H)$,
and $x_{k+1},\ldots,x_{\rho}$ such that $H=\Z v_P^{ }x_{k+1}\oplus\ldots
\oplus\Z v_P^{ }x_{\rho}\,$. Then $v_P^{ }K= \Z v_P^{ }x_1\oplus\ldots
\oplus\Z v_P^{ }x_{\rho}\,$. So if we choose the $y$'s as in
Theorem~\ref{hrwtd}, then we obtain that $F\subset K(x_1,\ldots,
x_{\rho},y_1,\ldots,y_{\tau})^h$. We set $K':=K(x_{k+1},\ldots,x_{\rho},
y_1,\ldots,y_{\tau})\sep$. Then it follows that $F':=F.K'\subset K(x_1,
\ldots,x_{\rho},y_1,\ldots,y_{\tau})^h.K'=K'(x_1,\ldots,x_k)^h$. (Here,
the last equality is seen as follows: $K'(x_1,\ldots,x_k)^h$ contains
$K'$, $K(x_1,\ldots,x_{\rho},y_1,\ldots,y_{\tau})$ and, by the universal
property of henselizations, also its henselization; hence,
``$\subseteq$'' holds. The converse follows from the universal property
since the left hand side is henselian, being an algebraic extension of
a henselian field.)

We can extend $Q$ to $K'(x_1,\ldots,x_k)^h$ in such a way that it
remains the identity on $K'$. By Lemma~\ref{prelBour},
\[K'(x_1,\ldots,x_k)Q\>=\>K'\;,\]
As $K'$ is separable-algebraically closed, $(K',\ovl{Q})$ is henselian.
Therefore, we can deduce from Lemma~\ref{hensPQ} and the fact that the
henselization is an immediate extension that $K'(x_1,\ldots,x_k)^hQ=K'$.
Since $K'\subseteq F'\subseteq K'(x_1,\ldots,x_k)^h$, it follows that
$F'Q= K'$. Hence by Corollary~\ref{corPQ}, $(F'|K'(x_1,\ldots,x_k),P)$
is uniformizable. By Proposition~\ref{prep}, the same holds for
$(K'(x_1,\ldots,x_k)|K',P)$. By transitivity, $(F'|K',P)$ is
uniformizable.

Pick $\zeta_1,\ldots,\zeta_m\in {\cal O}_F\,$. By what we have proved,
$(F'|K',P)$ is uniformizable with respect to $\zeta_1,\ldots,\zeta_m\,$.
By part c) of Lemma~\ref{going}, there is a finite Galois extension
${\cal K}'$ of $K(x_{k+1},\ldots,x_{\rho},y_1,\ldots, y_{\tau})$ such
that $(F.{\cal K}'|{\cal K}',P)$ is uniformizable with respect to the
$\zeta$'s. We take $\zeta'_1,\ldots,\zeta'_{m'}\in {\cal O}_{{\cal K}'}$
to be the uniformization coefficients of $F.{\cal K}',\zeta_1,\ldots,
\zeta_m$ in ${\cal K}'$. Since $P$ coincides with $\ovl{Q}$ on
${\cal K}'$, we know that $P$ is an Abhyankar place of ${\cal K}'|K$
of rank $r-1$. Now we have to distinguish two cases.

Suppose first that $({\cal K}',P)$ has rank $>1$. Then by induction
hypothesis, there is some ${\cal K}''$ obtained by at most $r-2$
(resp.\ $r-1$) Galois extensions from ${\cal K}'$, such that
$({\cal K}''|K,P)$ is uniformizable with respect to $\zeta'_1,\ldots,
\zeta'_{m'}\,$. We set ${\cal F}:=F.{\cal K}''= F.{\cal K}'.{\cal K}''$.
Then by part a) of Lemma~\ref{going} and the transitivity, $({\cal F}|
K,P)$ is uniformizable with respect to the $\zeta$'s. Further,
${\cal F}$ is obtained from $F$ by at most $r-1$ (resp.\ $r$)
Galois extensions, as required.

Now suppose that $({\cal K}',P)$ has rank $1$, i.e., $r=2$. By
Theorem~\ref{MTr1A} there is a finite Galois extension ${\cal K}''$ of
${\cal K}'$ such that $({\cal K}''|K,P)$ is uniformizable with respect
to $\zeta'_1,\ldots, \zeta'_{m'}\,$. As before, we set ${\cal F}:=
F.{\cal K}''$ and it follows that $({\cal F}|K,P)$ is uniformizable with
respect to the $\zeta$'s. Now ${\cal F}$ is obtained from $F$ by two
Galois extensions, as required.

If in addition $P$ is zero-dimensional, then we know from
Theorem~\ref{MTr1A} that there is a finite Galois extension $K'$ of
$K$ such that we can take ${\cal K}''={\cal K}'.K'$. In this case,
${\cal F}=F.{\cal K}''=F.{\cal K}'.K'$ is a Galois extension of $F$.
That is, ${\cal F}$ is obtained from $F$ by one Galois extension.

\bn
$\bullet$ \ {\bf Proof of Theorem~\ref{MTdens}}
\sn
Take function fields $F|K$ and $F'|K$ and a place $P$ of $F'|K$
such that $(F',P)$ lies in the completion of $(F,P)$. Then
by Proposition~\ref{comp1}, $(F'|F,P)$ is uniformizable.
Hence, the assertion of Theorem~\ref{MTdens} follows from the
corresponding assertions of Theorem~\ref{MTr1A} and Theorem~\ref{MTr>AG}
by use of part a) of Lemma~\ref{going} and transitivity.

\bn
\bn
\bn
{\bf References}
\newenvironment{reference}%
{\begin{list}{}{\setlength{\labelwidth}{5em}\setlength{\labelsep}{0em}%
\setlength{\leftmargin}{5em}\setlength{\itemsep}{-1pt}%
\setlength{\baselineskip}{3pt}}}%
{\end{list}}
\newcommand{\lit}[1]{\item[{#1}\hfill]}
\begin{reference}
\lit{[B]} {Bourbaki, N.$\,$: {\it Commutative algebra}, Paris (1972)}
\lit{[dJ]} {de Jong, A.~J.$\,$: {\it Smoothness, semi-stability and
alterations}, preprint}
\lit{[EL]} {Elliott, G.~A.$\,$: {\it On totally ordered groups, and
$K_0$}, in: Ring Theory Waterloo 1978, eds.\ D.~Handelman and
J.~Lawrence, Lecture Notes Math.\ {\bf 734}, 1--49}
\lit{[EN]} {Endler, O.$\,$: {\it Valuation theory}, Berlin (1972)}
\lit{[J--R]} {Jarden, M.\ -- Roquette, P.$\,$: {\it The Nullstellensatz
over $\wp$--adically closed fields}, J.\ Math.\ Soc.\ Japan {\bf 32}
(1980), 425--460}
\lit{[KA]} {Kaplansky, I.$\,$: {\it Maximal fields with valuations I},
Duke Math.\ J.\ {\bf 9} (1942), 303--321}
\lit{[KH--K]} {Khanduja, S.~K.\ -- Kuhlmann, F.--V.$\,$: {\it Valuations
on $K(x)$}, preprint, Toronto/ Sas\-ka\-toon (1997)}
\lit{[K1]} {Kuhlmann, F.--V.$\,$: {\it Henselian function fields and tame
fields}, preprint (extended version of Ph.D.\ thesis), Heidelberg
(1990)}
\lit{[K2]} {Kuhlmann, F.--V.$\,$: {\it Valuation theory of fields, abelian
groups and modules}, preprint, Heidelberg (1996), to appear in the
``Algebra, Logic and Applications'' series (Gordon and Breach), eds.\
A.~Mac\-intyre and R.~G\"obel}
\lit{[K3]} {Kuhlmann, F.--V.$\,$: {\it On places of algebraic function
fields in arbitrary characteristic}, submitted}
\lit{[K4]} {Kuhlmann, F.--V.$\,$: {\it Elementary properties of power
series fields over finite fields}, submitted; prepublication in:
Structures Alg\'ebriques Ordonn\'ees, S\'emi\-naire Paris VII (1997)}
\lit{[K5]} {Kuhlmann, F.--V.$\,$: {\it On local uniformization in
arbitrary characteristic}, The Fields Institute Preprint Series, Toronto
(1997)}
\lit{[K6]} {Kuhlmann, F.--V.$\,$: {\it On local uniformization in
arbitrary characteristic II}, in preparation}
\lit{[K7]} {Kuhlmann, F.--V.$\,$: {\it Valuation theoretic and model
theoretic aspects of local uniformization}, to appear in the
Proceedings of the Blowup Tirol Conference 1997}
\lit{[L]} {Lang, S.$\,$: {\it Algebra}, New York (1965)}
\lit{[R]} {Ribenboim, P.$\,$: {\it Th\'eorie des valuations}, Les
Presses de l'Uni\-versit\'e de Mont\-r\'eal (1964)}
\lit{[S]} {Spivakovsky, M.$\,$: {\it Resolution of singularities I:
local uniformization}, manu\-script, Toronto (1996)}
\lit{[W]} {Warner, S.$\,$: {\it Topological fields}, Mathematics
studies {\bf 157}, North Holland, Amsterdam (1989)}
\lit{[Z1]} {Zariski, O.$\,$: {\it Local uniformization on
algebraic varieties}, Ann.\ Math.\ {\bf 41} (1940), 852--896}
\lit{[Z2]} {Zariski, O.$\,$: {\it The reduction of singularities of an
algebraic surface}, Ann.\ Math.\ {\bf 40} (1939), 639--689}
\lit{[Z3]} {Zariski, O.$\,$: {\it A simplified proof for resolution of
singularities of an algebraic surface}, Ann.\ Math.\ {\bf 43} (1942),
583--593}
\lit{[Z--S]} {Zariski, O.\ -- Samuel, P.$\,$: {\it Commutative
Algebra}, Vol.\ II, New York--Heidel\-berg--Berlin (1960)}
\end{reference}
\adresse
\end{document}